\documentclass[12pt]{article}
\usepackage{amsfonts,amsmath}
\usepackage{amssymb,latexsym}
\usepackage[dvips]{epsfig}
\usepackage{amscd}
 
\title{Matrix factorizations and link homology II}
\author{ Mikhail Khovanov and Lev Rozansky}
\date{May 4, 2005}
 
\newtheorem{prop}{Proposition}
\newtheorem{theorem}{Theorem}
\newtheorem{lemma}{Lemma}

\newcommand{\oplusop}[1]{{\mathop{\oplus}\limits_{#1}}}
 
\begin{document} 

\maketitle
\baselineskip 14pt
 
\def\R{\mathbb R}
\def\Q{\mathbb Q}
\def\Z{\mathbb Z}
\def\l{\lbrace}
\def\r{\rbrace}
\def\o{\otimes}
\def\lra{\longrightarrow}
\def\Hom{\mathrm{Hom}}
\def\RHom{\mathrm{RHom}}
\def\Id{\mathrm{Id}}
\newcommand{\define}{\stackrel{\mbox{\scriptsize{def}}}{=}}
\def\drawing#1{\begin{center}\epsfig{file=#1}\end{center}}
 
\begin{abstract} To a presentation of an oriented link as the closure of a 
braid 
we assign a complex of bigraded vector spaces. The Euler characteristic of 
this complex (and of its triply-graded cohomology groups) 
is the HOMFLYPT polynomial of the link. We show that the dimension of 
each cohomology group is a link invariant. 
\end{abstract}

  
\section{Matrix factorizations with a parameter}\label{section-one}

In the paper [KR] we constructed, for each $n>0,$ 
a bigraded cohomology theory of links in $\R^3$ whose Euler characteristic 
is a certain one-variable specialization $(q^n,q)$ of the HOMFLYPT polynomial
[HOMFLY], [PT]. The $n=0$ specialization is the Alexander polynomial, 
equal to the Euler characteristic of the knot homology theory discovered by 
Ozsv\'ath, Rasmussen and Szab\'o [OS], [R1]. The approach in [KR] fails for $n=0,$ 
assigning trivial groups to any link. 

In this sequel to [KR] we assume that the reader is familiar 
with that paper. Recall that our construction of link cohomology was based 
on matrix factorizations with 
potentials being sums and differences of $x^{n+1},$ for various $x.$ 
When $n=0,$ the category of matrix factorizations (up to chain homotopies) 
with the potential $\sum \pm x_i$ is trivial. Looking for a remedy, let us add 
a formal variable $a$ and change the potential from $x$ to $ax.$ 

\begin{figure} \label{arc} \drawing{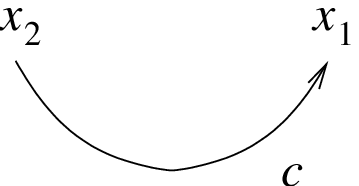} \caption{An arc} 
\end{figure}

Take an oriented arc $c$ as in figure~\ref{arc}, 
 label its ends $x_1$ and $x_2,$ 
and assign the potential $ax_1 - ax_2$ to the arc. Let 
$R=\Q[a,x_1,x_2]$ and define $C_c$ as the factorization
 $$ R \stackrel{a}{\lra} R \xrightarrow{x_1 -x_2} R.$$ 
We have $d^2= ax_1 - ax_2$ and view $C_c$ as an object 
of the homotopy category of matrix factorizations with the 
potential $a(x_1-x_2).$ 

Make $R$ bigraded by setting 
\begin{equation}\label{degrees}\deg(a) = (2,0), \hspace{0.3in} \deg(x_i)=(0,2).
\end{equation} 
This implies $\deg(d^2) = (2,2)$ and we select the bigrading 
of the middle $R$ in the factorization so that 
$\deg(d)= (1,1)$: 
 \begin{equation}\label{arcfact} 
  R \stackrel{a}{\lra} R\{-1,1\} \stackrel{x_1 -x_2}{\lra} R,\end{equation}  
where the bidegree shift by $(n_1,n_2)$ is denoted $\{n_1,n_2\}.$ 

\begin{figure} \drawing{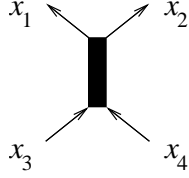} 
\caption{Wide edge $t$} \label{xonefour} \end{figure}

Next, given a wide edge $t$ as in figure~\ref{xonefour}, 
assign variables 
$x_1,x_2,x_3,x_4$ to the edges next to it. We can write 
$$ ax_1 + ax_2 - ax_3 -ax_4 = a(x_1+x_2-x_3-x_4) + 0 (x_1 x_2-x_3x_4).$$ 
Define $C_t$ as the tensor product (over $R$) of factorizations 
\begin{equation}
 R \stackrel{a}{\lra} R\{-1,1\} \xrightarrow{x_1+x_2-x_3-x_4} R\end{equation} 
and 
\begin{equation} 
 R\stackrel{0}{\lra} R\{-1,3\} \xrightarrow{x_1 x_2 - x_3 x_4} R\end{equation}  
where $R= \Q[a,x_1,x_2,x_3,x_4].$ 

Throughout the paper we work with matrix factorizations with 
potentials $w=a \sum_i \epsilon_i x_i$ where $i$ ranges over some finite set $I$ 
of integers and $\epsilon_i\in \{1,-1\}$ are "orientations" of $x_i.$ The category 
$mf_w$ has objects $(M,d)$ where $M=M^0\oplus M^1$ and $M^0,M^1$ 
are free bigraded $R$-modules (possibly of infinite rank), while $d$ is a generalized 
differential 
$$ M^0 \stackrel{d}{\lra} M^1 \stackrel{d}{\lra} M^0$$ 
of bidegree $(1,1)$ and subject to $d^2=w.$ Here $R$ is the ring of 
polynomials in $a$ and $x_i$'s with rational coefficients. The bidegrees are 
given by formula (\ref{degrees}).  Morphisms in $mf_w$ are bidegree-preserving 
maps of $R$-modules $M^0\to N^0,$ $M^1\to N^1$ that commute with 
$d.$ 

\begin{figure} \drawing{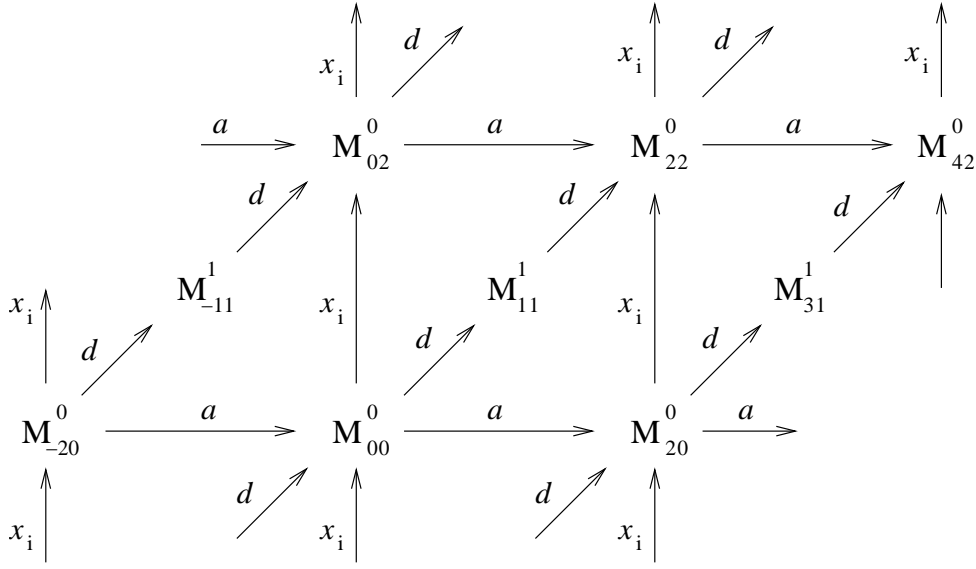} \caption{Lattice of factorization $M$} \label{lattice}
\end{figure}

We found it useful to visualize a matrix factorization as above by 
decomposing 
$$M^0=\oplusop{k,l} M^0_{k,l}, \hspace{0.2in} 
 M^1=\oplusop{k,l} M^1_{k,l},$$ 
as direct sums of vector spaces, 
one for each bidegree $(k,l),$ and placing 
them in the nodes of a coordinate lattice, see figure~\ref{lattice}. 
Diagonal arrows denote the differential, horizontal arrows show multiplication 
by $a$ and vertical arrows--multiplications by $x_i.$ 

The category $hmf_w$ of matrix factorizations up to chain homotopies has 
the same objects as $mf_w$ and the $\Q$-vector space of morphisms from 
$M$ to $N$ is the quotient of the space of morphisms in $mf_w$ by null-homotopic 
morphisms (the homotopy maps must have bidegree $(-1,-1)$). 

If the index set $I$ is empty, then $R=\Q[a]$ and $mf_w$ is equivalent to 
the category of complexes of free graded $\Q[a]$-modules. 

\begin{figure} \drawing{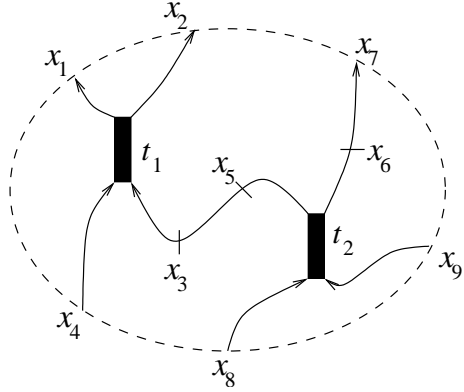} 
\caption{Graph $\Gamma,$ an example}\label{diagbound} \end{figure}

In general, given a planar marked graph $\Gamma$ (as described in [KR, Introduction]), 
possibly with boundary points, we assign to it a  matrix factorization $C(\Gamma)$ which 
is the tensor product of $C_c,$ over all arcs $c$ in $\Gamma,$ and $C_t,$ over all 
wide edges $t$ in $\Gamma.$ For instance, for the graph in 
figure~\ref{diagbound}, 
$$C(\Gamma) = C_{t_1} \otimes C_{t_2} \otimes C_{c_1} \otimes C_{c_2}$$ 
where $c_1,c_2$ are the arcs of $\Gamma$ with endpoints labelled 
$x_3,x_5$ and 
$x_7,x_6,$ respectively. The tensor product is taken over suitable polynomial 
rings $\Q[a,x_i]$ so that $C(\Gamma)$ is a finite rank free 
$\Q[a,x_1, \dots, x_9]$-module. The potential 
$$w=a(x_1 + x_2 - x_7 - x_4 - x_8 -x_9)$$ 
and we view $C(\Gamma)$ as an object of $mf_w$ (or $hmf_w$) with the 
ground ring $R$ 
 the polynomial ring $\Q[a,x_1,x_2,x_4,x_7,x_8,x_9]$ in $a$ and 
external (or boundary) variables. The other variables 
$x_3,x_5,x_6$ are "internal". Notice that $C(\Gamma)$ has 
infinite rank as an $R$-module.  

When $\Gamma$ has no boundary points, $w=0$ and $C(\Gamma)$ becomes 
a 2-periodic complex 
 $$ C^0(\Gamma) \stackrel{d}{\lra} C^1(\Gamma) \stackrel{d}{\lra} C^0(\Gamma)$$ 
of bigraded $\Q[a]$-modules. Its cohomology, denoted $H(\Gamma),$ is a 
bigraded $\Q[a]$-module.  

If $\Gamma$ is a single circle with one mark (glue together 
the endpoints of the arc in the figure~\ref{arc} and place a mark there), the complex 
is 
$$ \Q[a, x_1] \stackrel{a}{\lra}\Q[a,x_1]\{-1,1\} \stackrel{0}{\lra} \Q[a,x_1] $$ 
(since now $x_2=x_1$), and  $H(\Gamma)\cong \Q[x_1]\{-1,1\}.$ 

\begin{figure}  \drawing{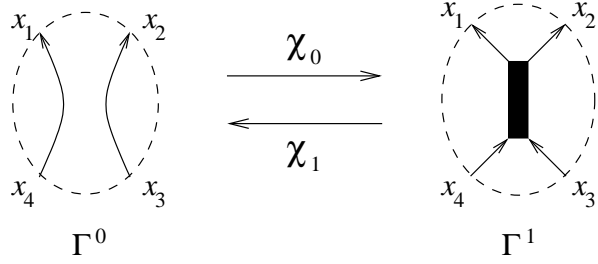}\caption{Graphs $\Gamma^0$
and $\Gamma^1$} \label{pair2}
 \end{figure}

Consider the diagrams $\Gamma^0, \Gamma^1$ in figure~\ref{pair2}. 
Factorization $C(\Gamma^0)$ is the tensor product of 
$$ R \stackrel{a}{\lra} R\{-1,1\} \xrightarrow{x_1-x_4} R $$
and 
$$ R \stackrel{a}{\lra} R\{-1,1\} \xrightarrow{x_2-x_3} R,$$
where $R=\Q[a,x_1,x_2,x_3,x_4].$ In the product basis, 
$C(\Gamma^0)$ has the form 
\begin{equation*}
  \begin{array}{l}R \\  \oplus \\ R\{-2,2\}\end{array}
  \stackrel{P_0}{\lra}
  \begin{array}{l} R\{-1,1\} \\ \oplus \\ R\{-1,1\}  \end{array}
  \stackrel{P_1}{\lra} 
  \begin{array}{l}R \\ \oplus \\ R\{-2,2\}\end{array}
\end{equation*}
with 
\begin{equation*}
P_0=\left(\begin{array}{cc} a  & x_3-x_2 \\
      a   & x_1-x_4 \end{array} \right), \hspace{0.2in}
P_1=\left(\begin{array}{cc} x_1-x_4  & x_2-x_3 \\
      -a  & a \end{array} \right).
\end{equation*}

Likewise,  $C(\Gamma^1)$ has the presentation 
\begin{equation*}
  \begin{array}{l}R \\ \oplus \\ R\{-2,4\}\end{array}
  \stackrel{Q_0}{\lra}
  \begin{array}{l} R\{-1,1\} \\ \oplus \\ R\{-1,3\}  \end{array}
  \stackrel{Q_1}{\lra} 
  \begin{array}{l}R \\ \oplus \\ R\{-2,4\}\end{array}
\end{equation*}
with 
\begin{equation*}
Q_0=\left(\begin{array}{cc} a  & x_3x_4-x_1x_2 \\
      0   & x_1+x_2-x_3-x_4 \end{array} \right), \hspace{0.2in}
Q_1=\left(\begin{array}{cc} x_1+x_2-x_3-x_4  & x_1x_2-x_3x_4 \\
      0  & a \end{array} \right). 
\end{equation*}
A map between $C(\Gamma^0)$ and $C(\Gamma^1)$ can be 
described by a pair of $2\times 2$ matrices with coefficients in $R$ 
that specify the images of the basis vectors of $C^i(\Gamma^0)$ 
in $C^i(\Gamma^1)$ for $i=0,1.$ 

Let  $\chi_0: C(\Gamma^0) \lra C(\Gamma^1)$ be 
given by the pair of matrices 
\begin{equation}\label{chinull} 
U_0^0=\left(\begin{array}{cc} x_4-x_2  &  0  \\
      0   &  1  \end{array} \right), \hspace{0.2in}
U_0^1=\left(\begin{array}{cc} x_4  &  -x_2 \\
      -1  & 1 \end{array} \right). 
\end{equation}
Our bases in $C(\Gamma^0)$ and $C(\Gamma^1)$ are homogeneous 
with respect to the bigrading of $R.$ It's easy to see that $\chi_0$ is 
a homogeneous map of bidegree $(0,2).$ 

Next, define $\chi_1: C(\Gamma^1) \lra C(\Gamma^0)$
by the pair of matrices 
\begin{equation}\label{chione} 
U_1^0=\left(\begin{array}{cc} 1  &  0  \\
      0   &  x_4-x_2  \end{array} \right), \hspace{0.2in}
U_1^1=\left(\begin{array}{cc} 1  &  x_2 \\
      1  & x_4 \end{array} \right). 
\end{equation}
The map $\chi_1$ is bidegree-preserving. 

\vspace{0.1in} 

\begin{figure} \drawing{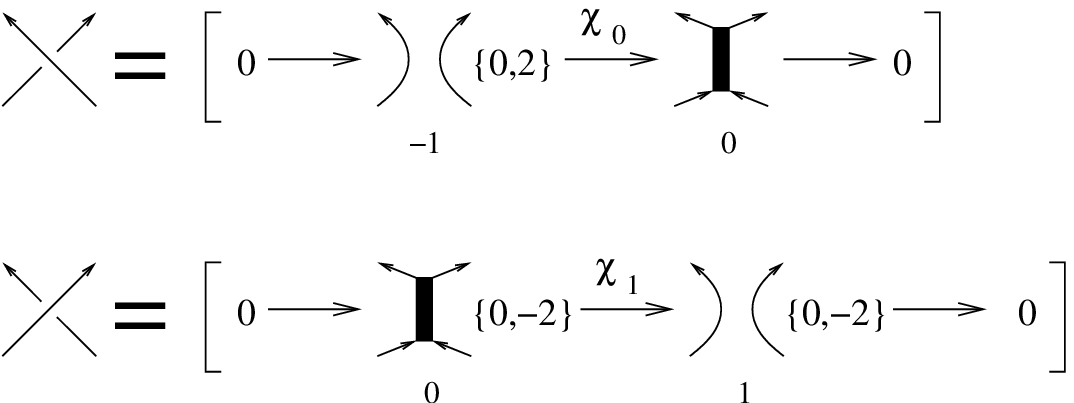}\caption{Complex assigned
 to a crossing} \label{cones}
 \end{figure}

Given a plane diagram $D$ of a tangle, place at least one mark on each 
internal edge of the diagram (an edge disjoint from the boundary of $D$), 
and label the marks and boundary points by $x_1, \dots, x_m.$ 
To each crossing $p$ of the diagram assign the 
 complex $C_p$ of matrix factorizations as follows. Up to shifts, 
the complex is the cone of the map $\chi_0$ or $\chi_1,$ depending 
on whether the crossing is 
positive or negative. The shifts are explained in figure~\ref{cones}. 

Thus, if the crossing is positive, 
$$C_p = \hspace{0.3in} 0 \lra C(\Gamma^0)\{0,2\} \stackrel{\chi_0}{\lra} 
 C(\Gamma^1) \lra 0 ,$$ 
with $C(\Gamma^1)$ positioned in cohomological degree $0.$ The shift $\{0,2\}$ 
makes the differential preserve the bidegree. 
 If the crossing is negative, 
$$ C_p = \hspace{0.3in} 0 \lra C(\Gamma^0)\{0,-2\} \stackrel{\chi_0}{\lra} 
 C(\Gamma^1)\{0,-2\} \lra 0 ,$$ 
with $C(\Gamma^0)$ in cohomological degree $0.$ The overall bigrading 
shift by $\{0,-2\}$ is here for the normalization of the Reidemeister 
move IIa (see later). 

Define $C(D)$ as the tensor product of $C_p,$ over all crossings $p$ of $D,$ 
and $C_c,$ over all arcs $c.$ It's a complex built out of matrix 
factorizations $C(\Gamma),$ over all resolutions $\Gamma$ of $D.$ 
The differential $\partial$ preserves the bigrading of each 
term $C^j(D).$ We view $C(D)$ as an object of the category $K(hmf_w).$  
The latter is the category whose objects are complexes of objects 
in $hmf_w$ and whose morphisms are homomorphism of complexes modulo null-homotopic 
morphisms.  

\vspace{0.1in} 

Now we specialize to the case when $D$ is a link diagram (has empty boundary). 
Then each term $C^j(D)$ in the complex $C(D)$ is an object of 
the homotopy category of bigraded free $\Q[a]$-modules.
We'll  see that $C^j(D),$ for any diagram $D,$ decomposes  
as a direct sum of contractible pieces 
$$0 \lra \Q[a] \stackrel{1}{\lra} \Q[a] \lra 0 $$ 
and the cohomology $H(C^j(D)),$ which we denote $CH^j(D).$ 
 Moreover, $a$ acts trivially on $CH^j(D),$ so we can 
ignore the $\Q[a]$-module 
structure and think of it as a bigraded $\Q$-vector space, 
$$ CH^j(D) = \oplusop{k,l} CH^j_{k,l}(D).$$
The bigrading descends from the bigrading on matrix factorizations 
$C(\Gamma).$ 

Thus, to $D$ we assign the complex $CH(D)$ of bigraded $\Q$-vector spaces 
$$  \dots \stackrel{\partial}{\lra}  
  CH^j(D)  \stackrel{\partial}{\lra} CH^{j+1}(D) \stackrel{\partial}{\lra} \dots $$
As a $\Q$-vector space, $CH(D)$ is the direct sum of cohomology groups $H(\Gamma)$ 
of complexes $C(\Gamma),$ over all resolutions of $D.$ 

The cohomology $H(D) = H(CH(D),\partial)$ of the above complex is triply-graded, 
$$H(D) = \oplusop{j,k,l} H^j_{k,l}(D).$$ 

Of course, for the whole construction to be interesting, 
$H(D)$ should not depend on the choice of $D,$ given $L.$ 

\begin{prop}\label{prop-mark} Let $D$ be a marked tangle diagram. 
Then $C(D)$, as an object of $K(hmf_w),$ does not depend 
on the number of markings on each edge of $D.$ 
\end{prop} 

Proposition~\ref{prop-mark} is proved in the next section.

\vspace{0.1in} 

\begin{figure}  \drawing{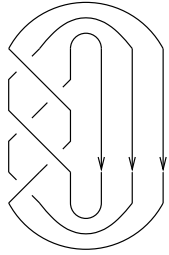}\caption{A braid diagram} 
 \label{braid} \end{figure}

Next we run into an obstacle: things seem to work well only if we 
restrict to diagrams $D$ that come from braids. Let's say that $D$ is 
a braid diagram 
if $D$ depicts the link $L$ as the closure of a clockwise oriented 
braid, see figure~\ref{braid}. We denote the braid by $D$ as well. 

\begin{figure}  \drawing{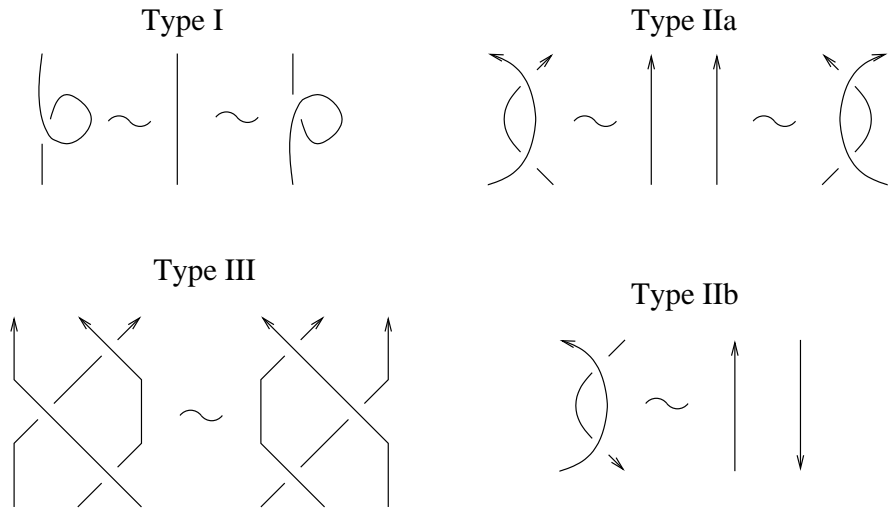}\caption{Reidemeister moves} \label{reid}
 \end{figure}

To justify the introduction of braid diagrams, partition the  
Reidemeister moves of links into types I, II, III in the usual way and 
then separate II into two subtypes, IIa and IIb, depending on 
orientations, see figure~\ref{reid}. We only consider 
the type III move with the orientations pointing in the same direction.  

\begin{prop} If diagrams $D_1, D_2$ are related by a Reidemeister move 
of type I, IIa or III, the complexes $C(D_1)$ and $C(D_2)$ are isomorphic 
as objects of the category $K(hmf_w),$ up to an overall shift in 
the triple grading.   
\end{prop} 

For a proof see Section~\ref{sec-proofs}. 

\vspace{0.1in} 

It seems likely to us that $C(D_1)$ and $C(D_2)$ are not isomorphic
(up to a shift)  
if $D_1, D_2$ are the two tangles in the IIb move in figure~\ref{reid}, 
although we did not prove this. To avoid the IIb move, we restrict our 
consideration to braid diagrams. Closures of braid diagrams $D_1, D_2$ 
are isotopic as oriented links iff $D_1, D_2$ are 
related by a chain of Markov moves. A Markov move is one of the following:  

(a) conjugation $DD' \leftrightarrow D'D,$

(b) transformations in the braid group:  
\begin{eqnarray*} 
D\sigma_j \sigma_i & \leftrightarrow & D \sigma_i \sigma_j \hspace{0.1in} \mathrm{if} 
 \hspace{0.1in} |i-j|>1, \\ 
D & \leftrightarrow & D \sigma_i \sigma_i^{-1},  \\
D & \leftrightarrow & D \sigma_i^{-1} \sigma_i,  \\
D \sigma_i \sigma_{i+1}\sigma_i & \leftrightarrow & D \sigma_{i+1}
  \sigma_{i} \sigma_{i+1} .
\end{eqnarray*} 

(c) transformations $D \leftrightarrow D \sigma_n^{\pm 1},$ 
for a braid $D$ with $n$ strands. 

\vspace{0.1in} 

Notice that we never see the Reidemeister move IIb when dealing 
with braid diagrams. The propositions stated earlier imply: 

\begin{theorem} 
Given two braid diagrams $D_1,D_2$ of an oriented link $L,$ 
the cohomology groups $H(D_1)$ and $H(D_2)$ are isomorphic as 
triply-graded vector spaces, up to an overall shift in the grading. 
\end{theorem} 

\begin{figure}  \drawing{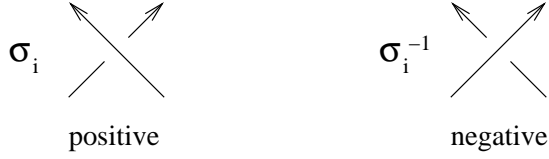}\caption{Positive and negative crossings 
and braid generators} \label{conv}
 \end{figure}

To describe the Euler characteristic of $H(D),$ we consider 
the function $F$ from braid diagrams to the ring of 
rational functions in $q$ and $t$ that is uniquely determined 
by the following properties: 

\begin{itemize} 
\item $F(D_1)=F(D_2)$ if $D_1,D_2$ are conjugate 
braid presentations (see Markov move (a) above),
\item $F(D_1)=F(D_2)$ if $D_1,D_2$ are related by a braid 
presentation move (Markov moves (b) above), 
\item $F(D \sigma_n) = F(D),$ for a braid $D$ with $n$ strands, 
\item $F(D\sigma_n^{-1})= - t^{-1}q^{-1} F(D),$ for a braid $D$ 
with $n$ strands, 
\item For any braid diagram $D$ there is a skein relation 
$$ q^{-1} F(D \sigma_i) - q F(D \sigma_i^{-1})= (q-q^{-1})
 F(D),$$ 
\item If $D$ is the one-strand diagram of the unknot, $F(D)=\frac{t^{-1}}{q^{-1}-q}.$ 
\end{itemize} 
To see that $F$ is simply a version of the 
HOMFLYPT polynomial, let $\alpha= - t^{-1}q^{-1}$ and 
consider
$$ 
\widetilde{F}(D) = \sqrt{\alpha}^{|D|_+ - |D|_-  - s(D)+1}  F(D) $$
where $|D|_+,$ respectively 
$|D|_-,$ is the number of positive, respectively negative,
crossings of $D,$ while $s(D)$ is the number of strands of $D.$ Our 
conventions are explained in figure~\ref{conv}. 
Then $\widetilde{F}(D)$ is invariant under all Markov moves of braids and 
satisfies the HOMFLYPT skein relation 
\begin{equation}\label{homfly-rel}
 q \sqrt{\alpha}\hspace{0.05in} \widetilde{F}(D\sigma_i^{-1}) - 
 (q \sqrt{\alpha})^{-1} \widetilde{F}(D\sigma_i) = (q-q^{-1}) \widetilde{F}(D).
\end{equation}  
Thus, $\widetilde{F}(D)$ equals the HOMFLYPT polynomial of the link $L,$ 
normalized so that 
$$\widetilde{F}(\mathrm{unknot})= \frac{\alpha }{1-q^{-2}}.$$

\begin{theorem} \label{sec-thm} 
For any braid diagram $D$ the Euler characteristic 
$\langle D \rangle$ of $H(D)$ equals $F(D).$
\end{theorem} 
The Euler characteristic 
\begin{equation*} \label{euler-char} 
  \langle D \rangle \define \sum_{j,k,l}(-1)^j t^k q^l \dim_{\Q} 
 H^j_{k,l}(D)
 \end{equation*} 
is a power series in $q$ with coefficients in $\Z[t,t^{-1}].$ 
The theorem claims that each vector space 
$H^j_{k,l}(D)$ is finite dimensional and the sum above 
is a rational function of $t$ and $q$ equal to $F(D).$ 
See the end of Section~\ref{sec-proofs} for a proof. 

\vspace{0.2in} 

Specializing to $q\sqrt{\alpha}=1$ in the equation 
(\ref{homfly-rel}) nets us the Alexander polynomial. In terms of 
$t$ and $q,$ we are imposing the relation $t=-q.$ Homologically, 
$t,$ $q$ and the minus sign correspond to the three grading 
directions. Hence, suitably collapsing the tri-grading to a
bigrading we get a categorification of the Alexander polynomial. 

\vspace{0.1in} 

Sergei Gukov, Albert Schwarz and Cumrun Vafa 
recently conjectured [GSV] 
that there exist integer-valued link invariants $D_{Q,s,r}$ 
depending on three integer parameters $Q,s,r,$ that can be 
used to determine ranks of  sl(n) link homology groups as 
well as the coefficients of the HOMFLYPT polynomial of a link. 
These invariants should come from the physical theory of the BPS 
states and should be related to ranks of cohomology 
groups of suitable moduli spaces. It would be interesting to try relating 
$D_{Q,s,r}$ to the dimensions of cohomology groups $H^j_{k,l}.$ 
Our normalization of the HOMFLYPT polynomial 
is similar to the one in [GSV], both having $q-q^{-1}$ as the 
denominator of the unknot invariant. 

On the other hand, it's been independently 
suggested by several people, including Oleg Viro [V], that there 
should exist a triply-graded link homology theory with the HOMFLYPT 
polynomial as the Euler characteristic. The current paper resulted 
from our search for such a theory and for a combinatorial 
categorification of the Alexander polynomial. 

Triply-graded cohomology theories had previously appeared in the 
work of Asaeda, Przytycki and Sikora [APS] on categorification of 
invariants of links in I-bundles over surfaces, and in 
 Audoux and Fiedler [AF], who introduced a refined Jones polynomial 
and its categorification, which are only invariant under braid-like 
isotopies. Restriction to  braid-like isotopies appears in our 
construction as well, but 
we don't know how our invariant relates to those of [APS] and [AF].  

\vspace{0.1in} 

To conclude this section, we mention several modifications,  
potential generalizations and illnesses of the homology theory $H.$ 

\begin{itemize} 
\item It's not natural that we have to restrict to braid diagrams
to get a link invariant. In another sign of disfunctionality, 
the theory does not extend to all cobordisms. For instance, the 
cohomology groups of the unknot do not have a Frobenius algebra 
structure over the cohomology ring of the empty link (it's 
convenient to define the latter ring to be $\Q[a]$), preventing us 
from extending the theory even to unknotted cobordisms between 
unlinks.  
\item Any field $k$ can be used instead of $\Q.$ More generally, we 
can work over $\Z,$ so that the invariant of a closed planar graph 
$\Gamma$ is a complex of graded free abelian groups, 
up to chain homotopy equivalence, and the invariant of a link is 
a complex of complexes as above, up to chain homotopy equivalence. 
Taking the homology $H(\Gamma,\Z)$ 
of each resolution of $D$ and forming a complex out of them produces 
a complex $CH(D,\Z)$ from a diagram $D.$ We then specialize to 
braid diagrams and take the cohomology of $CH(D,\Z).$ The resulting 
groups $H(D,\Z)$ are triply-graded and, up to isomorphism, do not 
depend on the choice of braid diagram $D,$ given $L.$ 
\item In Section~\ref{sec-proofs} we rewrite the factorizations 
$C(\Gamma^0), 
C(\Gamma^1)$ and the maps $\chi_0,$ $\chi_1$ in the form that depends only 
on $a$ and the differences $x_i-x_j$ of the variables $x_i.$ This allows 
us to pass from the ring $R=\Q[a,x_1, \dots, x_m]$ to the smaller ring 
$\overline{R}=\Q[a,x_2-x_1, \dots, x_m-x_1].$ The definition of cohomology 
and the proof of its invariance work over $\overline{R}$ as well, leading to 
\emph{reduced} cohomology groups $\overline{H}(D),$ with the property 
$H(D)= \overline{H}(D)\otimes_{\Q} \Q[x].$ In the reduced theory the 
unknot has one-dimensional cohomology groups. 
\item $sl(n)$ link homology theory (see [KR]) utilized the potential $x^{n+1}.$
 Soon afterwards Gornik [G] studied a deformation 
of that theory with the potential $x^{n+1}-(n+1)\beta^n x.$ In the $n=2$ case
the deformation was found earlier by Lee [L] and used by Rasmussen 
in his combinatorial proof of the Milnor conjecture [R2]. 
The definitions in [KR] can be generalized to the 
potential $x^{n+1}+a_n x^n + \dots + a_1 x$ where $a_1, \dots, a_n$ are 
formal variables. We hope that this generalization will be invariant under 
the Reidemeister moves and will turn out to be the "$sl(n)$-equivariant" version 
of $sl(n)$ link homology. The invariant of the empty link should be the 
ring of polynomials in $a_1, \dots, a_n,$ and naturally isomorphic to the 
$U(n)$-equivariant cohomology ring of the point. The invariant of the unknot 
should be the quotient of the polynomial ring $\Q[x,a_1, \dots, a_n]$ by 
the relation $x^{n+1}+a_n x^n + \dots + a_1 x=0,$ isomorphic to 
the $U(n)$-equivariant cohomology ring of $\mathbb{CP}^n.$ A certain 
version of Bar-Natan link homology [BN], [K]  should correspond to the  
 potential $x^3+ a x.$ 
\item For a common generalization of the $U(n)$-equivariant link homology 
and the theory described here one could try the potential 
$a_{n+1}x^{n+1}+a_n x^n + \dots + a_1 x$ with all $a$'s being 
formal variables. Factorizations $C(\Gamma^0), C(\Gamma^1),$ the maps 
$\chi_0, \chi_1$ and the complex $C(D)$ can be defined for this potential 
as well, but we don't know whether this theory will be invariant under 
the Reidemeister-Markov moves of braid diagrams. 
\end{itemize}


\section{Proofs} \label{sec-proofs}

\vspace{0.05in} 

{\bf 1. Product factorizations, graph homology and Koszul complexes.} 

Given a polynomial ring $R$ and a pair of elements $a_1,b_1\in R,$ we denote by 
$(a_1,b_1)$ the factorization 
$$ R \stackrel{a_1}{\lra} R \stackrel{b_1}{\lra} R.$$
Given a finite set of such pairs $(a_i,b_i), 1\le i\le n,$ we denote by 
$(\mathbf{a}, \mathbf{b})$ their tensor product (over $R$): 
 $$ ({\mathbf a}, {\mathbf b}) \define \otimes_{i} (a_i,b_i), \hspace{0.3in} 
 {\mathbf a}=(a_1, \dots, a_n), \hspace{0.1in} {\mathbf b}=(b_1, \dots, b_n).$$
We will also write $(\mathbf{a}, \mathbf{b})$ in the column form 
$$ \left( \begin{array}{cc} a_1 & b_1 \\
   a_2 & b_2 \\  \dots & \dots  \\ a_n & b_n \end{array} \right) $$
and call it a Koszul factorization. 
By an elementary transformation of rows $i$ and $j$ we mean a modification 
$$\left( \begin{array}{cc} a_i & b_i \\ a_j & b_j \end{array} \right) 
 \hspace{0.15in} \xrightarrow{[ij]_{\lambda}} \hspace{0.15in} 
 \left( \begin{array}{cc} a_i & b_i+\lambda b_j \\ a_j-\lambda 
 a_i & b_j \end{array} \right)$$
for some $\lambda\in R.$ We denote it by $[ij]_{\lambda}.$ 
All other rows of $(\mathbf{a},\mathbf{b})$ are left 
unchanged. An elementary transformation takes a Koszul 
factorization $(\mathbf{a},\mathbf{b})$
to an isomorphic factorization, since we're only changing a 
basis vector in the free $R$-module 
underlying the factorization $(\mathbf{a},\mathbf{b}).$

\vspace{0.1in} 

Suppose now that $y$ is one of the generators of the polynomial ring 
 $R,$ so we can write $R=R'[y],$ and that the potential   
 $w=\sum a_i b_i$ lies in $R'$ (in this situation we say that $y$ 
is an \emph{internal} variable). Then any 
 factorization $M$ over $R$ restricts to (an infinite rank) 
 factorization over $R',$ 
which we denote $M'.$ Assume furthermore that one of the rows in 
$(\mathbf{a},\mathbf{b})$ has the form $(0, y - \mu)$ where $\mu\in R'.$ 
Denote by $(\mathbf{a}',\mathbf{b}')$ the factorization over $R'$ obtained from 
$(\mathbf{a},\mathbf{b})$ by removing the row $(0, y - \mu)$ and substituting 
$\mu$ for $y$ everywhere in all other rows. 

\begin{prop} \label{without-row} Factorizations $(\mathbf{a}',\mathbf{b}')$ and 
$(\mathbf{a},\mathbf{b})'$ are chain homotopy equivalent. 
\end{prop} 

\emph{Proof:} By changing a variable $y\to y-\mu$ we reduce to the case $\mu=0.$ 
We can write  $a_i=a_i'+ y a_i''$ and $b_i=b_i'+y b_i''$ 
where $a_i', b_i'\in R'.$ 
Applying elementary transformations to rows $(0,y)$ and $(a_i,b_i)$ 
we reduce the latter to $(a_i,b_i'),$ while $(0,y)$ is transformed into 
$(\sum a_i b_i'', y).$ Next, 
change $(\sum a_i b_i'',y)$ into $(y,\sum a_i b_i'')$ (this shifts factorization 
 $M$ to $M\langle 1 \rangle$) and  
apply elementary transformations to rows $(a_i,b_i')$ and $(y,\sum a_i b_i'').$ 
The row $(a_i, b_i')$ becomes $(a_i', b_i'),$ while the row with $y$ turns
into $(y, \sum (a_i b_i'' + a_i'' b_i')).$ Since the potential does not depend on $y,$ 
the latter sum is zero. Now shift $(y,0)$ back to $(0,y).$ The result is the Koszul factorization, 
isomorphic to 
$(\mathbf{a},\mathbf{b}),$ with rows $(a_i',b_i')$ and $(0,y).$ 
This factorization 
is the tensor product of $(\mathbf{a}',\mathbf{b}'),$ as defined above, and 
$(0,y).$ Therefore, $(\mathbf{a},\mathbf{b})$ is isomorphic, as 
an $R$-factorization, to the total factorization of the bifactorization 
 $$ (\mathbf{a}',\mathbf{b}')\otimes_{R'}R'[y] \stackrel{0}{\lra}
 (\mathbf{a}',\mathbf{b}')\otimes_{R'}R'[y] \stackrel{y}{\lra} 
 (\mathbf{a}',\mathbf{b}')\otimes_{R'}R'[y]. $$ 
As a factorization over the smaller ring $R',$ it decomposes into a direct sum of 
 contractible factorizations which are the total factorizations of  
 $$  (\mathbf{a}',\mathbf{b}')\otimes  y^{j+1}  \stackrel{0}{\lra}
   (\mathbf{a}',\mathbf{b}')\otimes  y^{j}  \stackrel{y}{\lra} 
 (\mathbf{a}',\mathbf{b}')\otimes  y^{j+1}, $$
for $j\ge 0,$ and the factorization $ (\mathbf{a}',\mathbf{b}').$ Proposition 
follows. $\square$
 
\emph{Remark:} The second half of the above proof just says that the 
complex of $R'$-modules 
$$0\lra R'[y] \stackrel{y}{\lra} R'[y] \lra 0$$ 
is the direct sum of contractible complexes 
$$ 0\lra R' y^j \stackrel{y}{\lra} R' y^{j+1} \lra 0$$ 
and the complex $0 \lra R' \lra 0.$ 

\vspace{0.1in} 

Suppose we are given a planar marked graph $\Gamma,$ possibly with 
boundary. To $\Gamma$ we assigned a Koszul factorization $C(\Gamma)$
which has a rather special form. Each arc in $\Gamma$ contributes the  
row $(a , x_i-x_j)$ to the Koszul matrix of $C(\Gamma),$ where $x_i$ and 
$x_j$ are the labels at the endpoints of the arc. Each wide edge in $\Gamma$ 
contributes two rows 
$$\left( \begin{array}{cc}  a & x_i + x_j -x_k -x_l \\
                                           0  & x_i x_j - x_k x_l 
    \end{array} \right) $$
to the Koszul matrix, where $x_i, x_j,$ $x_k, x_l$ are the labels bounding 
the edge. If $\Gamma$ has $m_1$ arcs and $m_2$ wide edges, the Koszul 
matrix of $C(\Gamma)$ will have $n=m_1+2m_2$ rows. Permute these rows 
so that the first $m_1+m_2$ rows have the form $(a, z)$ where $z$'s are 
some linear functions of $x_i$'s. We call these rows \emph{linear} rows. 
The last $m_2$ rows have the form 
$(0 ,  x_i x_j - x_k x_l)$ for various quadruples of indices 
$(i,j,k,l).$ Call these \emph{quadratic} rows. 
 
Apply elementary transformations with $\lambda=1$ to the first row 
paired with every other linear row. In other words, we convert $b_1$ to 
$b_1+b_2+\dots + b_{m_1+m_2}$ and subtract $a_1=a$ from $a_p=a$ 
for $p=2,3,\dots, m_1+m_2.$ The Koszul matrix transforms 
into a matrix with the first row $(a, \sum \epsilon_i x_i)$ where 
the sum is over all boundary points of $\Gamma$ and $\epsilon_i = \pm 1$ 
depending on the orientation of $\Gamma$ at that point. All other linear 
rows acquire the form $(0, z),$ with the same linear functions $z$ as 
before. Nothing happens to the quadratic rows. The Koszul matrix now 
has the form 
$$\left( \begin{array}{cc} a  &  b_1 \\ 0 & b_2 \\ \cdots & \cdots \\
  0  & b_n \end{array} \right)$$
with $b_1= \sum \epsilon_i x_i.$
After this change of basis, every row but the first one has the 
first term zero. Hence, it comes from a one-term Koszul complex 
$$0 \lra R \stackrel{b_p}{\lra} R\lra 0$$ 
by collapsing cohomological grading from $\Z$ to $\Z_2.$ 
Likewise, the tensor product of all rows save the first is a factorization 
obtained from the Koszul complex of the sequence $(b_2, b_3,\dots , b_n)$ 
by collapsing the grading. 

Note that our polynomial ring is, in addition, bigraded. 
Taking all gradings into account, the collapse is from a 
triple grading to a bigrading (see figure~\ref{lattice}). 
No cyclic components appear in the collapsed 
grading since the differential has nonzero bidegree $(1,1).$ 
Finally, observe that in the new Koszul matrix parameter $a$ appears only once, 
in the first row. 

\vspace{0.1in} 

Next consider the case when $\Gamma$ is closed (has no boundary points). 
The first row becomes $(a\hspace{0.1in} 0)$ and the whole factorization 
comes from the Koszul complex of the sequence $(a, b_2, \dots, b_n).$ 
by collapsing its grading. Moreover, $a$ plays a purely decorative 
role, and, using proposition~\ref{without-row}, we can throw out this row 
simultaneously with removing $a$ from the list of variables, which would then 
have only $x_i$'s. In other words, the cohomology $H(\Gamma)$ of the 
factorization $C(\Gamma)$ is isomorphic to the cohomology of 
the Koszul complex of the sequence  $(b_2, b_3, \dots , b_n),$ 
with the trigrading collapsed to a bigrading. 

Thus, although the 2-periodic complex $C(\Gamma)$ 
as well as its cohomology $H(\Gamma)$ are $\Q[a]$-modules, $a$ acts 
trivially on $H(\Gamma).$ 

\vspace{0.1in}


{\bf 2. Maps $\chi_0, \chi_1$ revisited.} 

\vspace{0.1in} 

Recall the row operation $[ij]_{\lambda}$ on a Koszul matrix 
of a factorization: 
$$\left( \begin{array}{cc} a_i & b_i \\ a_j & b_j \end{array} \right) 
 \hspace{0.15in} \lra \hspace{0.15in} 
 \left( \begin{array}{cc} a_i & b_i+\lambda b_j \\ a_j-\lambda 
 a_i & b_j \end{array} \right)$$
Denote by $|0\rangle$ and $|1\rangle$ the standard basis vectors 
in factorizations 
$(a_i, b_i)$ and $(a_j, b_j)$: 
\begin{eqnarray*} & & R|0\rangle \stackrel{a_i}{\lra} R|1 \rangle 
 \stackrel{b_i}{\lra} R|0\rangle , \\
 & & R|0\rangle \stackrel{a_j}{\lra} R|1 \rangle 
 \stackrel{b_j}{\lra} R|0\rangle. 
\end{eqnarray*} 
Let $|00\rangle, |01\rangle, |10\rangle, |11\rangle$ be the standard basis 
vectors in the tensor product factorization $(a_i,b_i)\otimes (a_j,b_j).$ 
The row operation $[ij]_{\lambda}$ corresponds to the isomorphism of 
factorizations 
$$ (a_i,b_i)\otimes (a_j,b_j) \cong (a_i , b_i + \lambda b_j) \otimes 
  (a_j -\lambda a_i, b_j)$$ 
which takes the standard basis of the LHS factorization to the basis  
$$ |00\rangle, |01\rangle, |10> + \lambda |01\rangle, |11\rangle$$ 
of the RHS tensor product. 

Denote by $\psi(y)$ the following morphism between two Koszul factorizations: 
  \begin{equation*}
   \begin{CD}
     R @>{x}>> R @>{yz}>> R  \\
     @V{1}VV      @V{y}VV   @V{1}VV \\
     R @>{xy}>> R @>{z}>> R
   \end{CD}
 \end{equation*}

\begin{lemma}\label{comm-lemma}
  The following squares are commutative: 
  \begin{equation*}
   \begin{CD}
  \left( \begin{array}{cc} a_1 & b_1 \\ a_2 & b_2 c_2 \end{array}\right)   
   @>{\mathrm{Id}\otimes \psi(c_2)}>>  
 \left( \begin{array}{cc} a_1 & b_1 \\ a_2c_2 & b_2  \end{array}\right)   
   \\
     @V{[12]_{\lambda}}VV      @V{[12]_{\lambda c_2}}VV    \\
    \left( \begin{array}{cc} a_1 & b_1+\lambda b_2 c_2 \\ a_2-\lambda a_1 & 
    b_2 c_2 \end{array}\right)   
    @>{\mathrm{Id}\otimes \psi(c_2)}>> 
    \left( \begin{array}{cc} a_1 & b_1+\lambda b_2 c_2 \\ (a_2-\lambda a_1)c_2 & 
    b_2  \end{array}\right)   
   \end{CD}
 \end{equation*}
 \begin{equation*}
   \begin{CD}
  \left( \begin{array}{cc} a_1 & b_1 \\ a_2 & b_2 c_2 \end{array}\right)   
   @>{\mathrm{Id}\otimes \psi(c_2)}>>  
 \left( \begin{array}{cc} a_1 & b_1 \\ a_2c_2 & b_2  \end{array}\right)   
   \\
     @V{[21]_{\lambda c_2}}VV      @V{[21]_{\lambda}}VV    \\
    \left( \begin{array}{cc} a_1-\lambda c_2 a_2 & b_1 \\ 
   a_2 & 
    c_2 (b_2+\lambda b_1)  \end{array}\right)   
    @>{\mathrm{Id}\otimes \psi(c_2)}>> 
    \left( \begin{array}{cc} a_1-\lambda c_2 a_2 & b_1 \\ 
   a_2 c_2  & 
    b_2 +\lambda b_1 \end{array}\right)   
   \end{CD}
 \end{equation*}
\end{lemma} 

\emph{Proof:} direct computation. $\square$ 

\vspace{0.1in} 

Denote by $\psi'(y)$ the "opposite" morphism of $\psi(y)$: 
  \begin{equation*}
   \begin{CD}
     R @>{xy}>> R @>{z}>> R  \\
     @V{y}VV      @V{1}VV   @V{y}VV \\
     R @>{x}>> R @>{yz}>> R
   \end{CD}
 \end{equation*}
The analogue of lemma~\ref{comm-lemma} holds for $\psi'$ 
as well (just reverse all horizontal arrows in the commutative 
diagrams above).  We call $\psi$ and $\psi'$
\emph{flip} morphisms. 

\vspace{0.1in} 

Starting with the Koszul matrices for $C(\Gamma^0)$ and 
$C(\Gamma^1)$ and 
applying a row transformation to each of them, we get the 
following equivalent Koszul forms for these factorizations: 
\begin{eqnarray} 
C(\Gamma^0): &  \label{fla-gam-n}
\left( \begin{array}{cc} a & x_1-x_4 \\ a & x_2 -x_3 \end{array} 
   \right)\hspace{0.06in} \xrightarrow{[12]_1} \hspace{0.06in} 
  \left( \begin{array}{cc} a & x_1+x_2-x_3-x_4 \\ 0 & x_2 -x_3 \end{array} 
   \right) &   \\
 C(\Gamma^1): & \label{fla-gam-o}  
\left( \begin{array}{cc} a & x_1+x_2-x_3-x_4 \\ 0 & x_1x_2-x_3x_4 \end{array} 
   \right)\hspace{0.06in} \xrightarrow{[21]_{-x_2}} \hspace{0.06in} 
  \left( \begin{array}{cc} a & x_1+x_2-x_3-x_4 \\ 0 & (x_2 -x_3)(x_4-x_2) 
  \end{array} 
   \right) & 
\end{eqnarray} 
The first rows of these new Koszul matrices for $C(\Gamma^0), C(\Gamma^1)$ 
are identical while the second rows look related. In fact, there is a 
flip morphism $\psi(x_4-x_2)$ from $(0,(x_2 -x_3)(x_4-x_2))$ to 
$(0, x_2 -x_3)$: 
  \begin{equation*}
   \begin{CD}
     R @>{0}>> R @>{(x_2-x_3)(x_4-x_2)}>> R  \\
     @V{1}VV      @V{x_4-x_2}VV   @V{1}VV \\
     R @>{0}>> R @>{\quad  x_2-x_3 \quad }>> R
   \end{CD}
 \end{equation*}
and the flip morphism $\psi'(x_4-x_2)$ back. Tensoring these flip 
morphisms with the identity morphism on the first row, we obtain 
maps of factorizations 
$$ \mathrm{Id}\otimes \psi'(x_4-x_2)\hspace{0.02in} : \hspace{0.02in} 
  C(\Gamma^0)\lra C(\Gamma^1), \hspace{0.35in} 
   \mathrm{Id}\otimes \psi(x_4-x_2)\hspace{0.02in} : \hspace{0.02in} 
  C(\Gamma^1)\lra C(\Gamma^0). $$

\begin{lemma} Maps $\mathrm{Id}\otimes \psi'(x_4-x_2)$ and 
 $ \mathrm{Id}\otimes \psi(x_4-x_2)$ are equal to $\chi_0$ and $\chi_1,$ 
 respectively. 
\end{lemma} 
The proof is a straightforward linear algebra computation. 
$\square$ 

\vspace{0.1in} 

Therefore, our definition of the complex $C(D)$ of factorizations 
assigned to a tangle diagram can be rewritten via modified Koszul 
matrices as above and maps $\psi, \psi'.$ We'll use this alternative 
presentation in our proof of the invariance of $C(D)$ below. 
The new definition simplifies the appearance of $C(D)$ by creating 
more zeros in the Koszul matrices of $C(\Gamma)$ and making the
differential easier to describe and understand (at the cost of 
breaking the "lateral" symmetry $x_1\leftrightarrow x_2,$ 
$x_3 \leftrightarrow x_4$ of the original Koszul matrices). 
The differential acts now as the identity on all but $m_2$ rows, 
where $m_2$ is the number of crossings of $D.$  

\vspace{0.1in} 


{\bf 3. Markings don't matter.} 

\vspace{0.1in} 

To define the complex $C(D)$ for a tangle diagram $D,$ we need to 
place several marks on $D$: at least one on each internal edge and 
each circle and some (possibly none) on each external edge (an edge 
containing a boundary point). In this subsection we prove 
proposition~\ref{prop-mark} that was stated earlier and says that, 
up to chain homotopy equivalence, $C(D)$ does not depend on how 
marks are placed on the edges of $D.$ 

\begin{figure}  \drawing{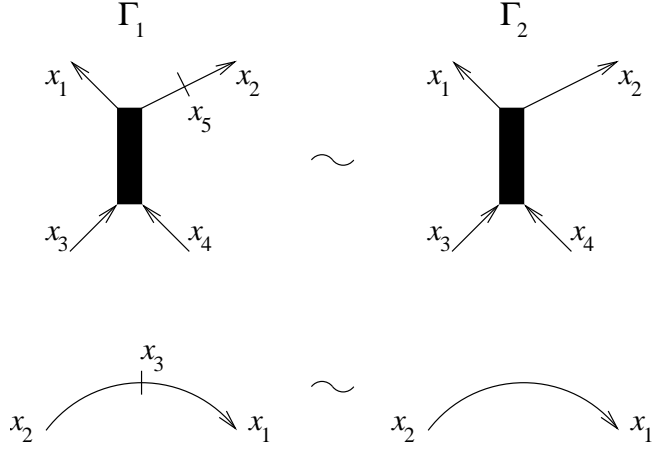}\caption{Mark removal equivalences} 
\label{rmarks1}
 \end{figure}

\begin{lemma} Factorizations $\Gamma_1$ and $\Gamma_2$ 
are isomorphic in $hmf_w$ if $\Gamma_2$ is obtained from $\Gamma_1$ 
by removing a mark. 
\end{lemma} 

\emph{Proof:} It's enough to check this property locally. We depicted 
two such local pairs $(\Gamma_1, \Gamma_2)$ in 
figure~\ref{rmarks1}, and refer the reader to [KR] for a more 
detailed treatment. We only check the isomorphism for the top 
pair in figure~\ref{rmarks1}, other cases are similar. We transform 
the Koszul matrix of $C(\Gamma_1)$ as follows:  
$$ \left( \begin{array}{cc} a & x_1+x_5-x_3-x_4 \\ 0 & 
  x_1 x_5 -x_3 x_4 \\  a & x_2 - x_5 \end{array} 
  \right) \stackrel{[13]_1}{\lra} 
  \left( \begin{array}{cc} a & x_1+x_2-x_3-x_4 \\ 0 & 
  x_1 x_5 -x_3 x_4 \\  0 & x_2 - x_5\end{array} \right)$$
The variable $x_5$ is internal. According to 
proposition~\ref{without-row} with $y=x_5$ we can remove the last 
row of the RHS matrix, substitute $x_2$ for $x_5$ everywhere else and 
forget about $x_5.$ We end up with the Koszul matrix of $C(\Gamma_2).$ 
$\square$ 

\begin{figure}  \drawing{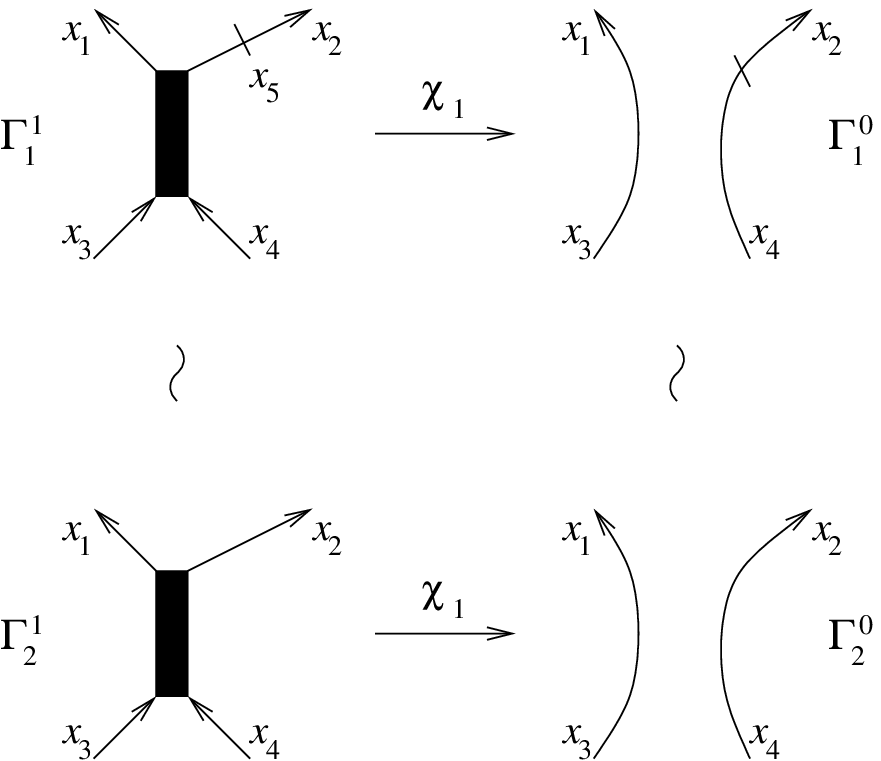}\caption{} \label{rmarks2}
 \end{figure}

\vspace{0.1in} 

To show independence of $D$ on the number and position of marks, we 
need to check compatibility of the isomorphisms $C(\Gamma_1)\cong C(\Gamma_2)$ 
above with maps $\chi_0, \chi_1.$ We'll only work through one case and 
leave the others to an interested reader. Let's check that complexes 
of factorizations 
$$0 \lra C(\Gamma^1_1) \stackrel{\chi_1}{\lra}  C(\Gamma^0_1)\lra 0 $$ 
and 
\begin{equation} 
0 \lra   C(\Gamma^1_2) \stackrel{\chi_1}{\lra}C(\Gamma^0_2)\lra 0
\label{fla-sec} 
\end{equation}  
are chain homotopy equivalent, for figure~\ref{rmarks2} diagrams. 
The first complex, written via Koszul matrices, has the form 
$$ \left( \begin{array}{cc} a & x_1 +x_5 -x_3 -x_4 \\
   0 & (x_5 - x_3)(x_4-x_5) \\ a & x_2 -x_5 \end{array}\right) 
 \xrightarrow{\Id \otimes \psi(x_4-x_5) \otimes Id} 
   \left( \begin{array}{cc} a & x_1 +x_5 -x_3 -x_4 \\
   0 & x_5-x_3 \\ a & x_2 -x_5 \end{array}\right) $$
Applying $[13]_1$ simultaneously to both matrices we get an 
isomorphic complex 
$$\left( \begin{array}{cc} a & x_1 +x_2 -x_3 -x_4 \\
   0 & (x_5 - x_3)(x_4-x_5) \\ 0 & x_2 -x_5 \end{array}\right) 
 \xrightarrow{\Id \otimes \psi(x_4-x_5) \otimes Id} 
   \left( \begin{array}{cc} a & x_1 +x_5 -x_3 -x_4 \\
   0 & x_5-x_3 \\ 0 & x_2 -x_5 \end{array}\right) $$
The only internal variable is $x_5.$ We switch from $x_5$ to 
$x=x_2-x_5.$ We think of $x$ as an internal variable, 
while $a, x_1, x_2,$ $x_3, x_4$ are external. Both matrices 
have identical bottom rows $(0,x)$ and the differential 
is the identity on that row. Therefore, we can eliminate $x$ from 
the complex, reducing the ground ring to $\Q[a,x_1, x_2, x_3, x_4],$ 
crossing out the bottom row and setting $x=0.$ The resulting 
complex  
$$\left( \begin{array}{cc} a & x_1 +x_2 -x_3 -x_4 \\
   0 & (x_2 - x_3)(x_4-x_2) \end{array}\right) 
 \xrightarrow{\Id \otimes \psi(x_4-x_2) } 
   \left( \begin{array}{cc} a & x_1 +x_2 -x_3 -x_4 \\
   0 & x_2-x_3 \end{array}\right) $$
is isomorphic to (\ref{fla-sec}). $\square$ 

\vspace{0.1in} 


{\bf 4. Invariance under Reidemeister move I.} 

\begin{figure}  \drawing{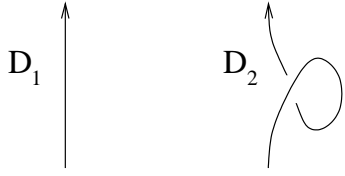}\caption{Type IA Reidemeister move} 
 \label{curl1} \end{figure}
\begin{figure} \drawing{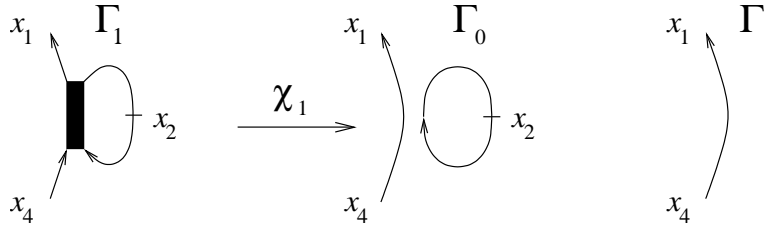} \caption{Resolutions of $D_2$} 
 \label{digon} \end{figure}

Consider type IA Reidemeister move as depicted in figure~\ref{curl1}. 
The complex $C(D_2)\{0,2\}$ has the form 
$$ 0 \lra C(\Gamma_1) \stackrel{\chi_1}{\lra} C(\Gamma_0) \lra 0, $$ 
see figure~\ref{digon}. In terms of Koszul matrices, the complex is 
given by  
$$ \left( \begin{array}{cc} a & x_1 -x_4 \\ 0 & 0 \end{array} 
 \right) \xrightarrow{\Id \otimes \psi(x_4-x_2)}
  \left( \begin{array}{cc} a & x_1 -x_4 \\ 0 & 0 \end{array} 
 \right) $$
(we set $x_3=x_2$ in the formulas (\ref{fla-gam-n}), (\ref{fla-gam-o}) 
for $\Gamma^0$ and $\Gamma^1$).

The differential is the identity on the first row, and on the second 
row given by 
 \begin{equation*}
   \begin{CD}
     R @>{0}>> R\{-1,3\} @>{0}>> R  \\
     @V{1}VV      @V{x_4-x_2}VV   @V{1}VV \\
     R @>{0}>> R\{-1,1\} @>{0}>> R
   \end{CD}
 \end{equation*}
Therefore, the complex splits into a direct sum of a contractible 
complex and the tensor product of $(a, x_1-x_4)$ with 
 $$ 0 \lra R\{-1,3\} \xrightarrow{x_4-x_2} R\{-1,1\} \lra 0.$$ 
Since $x_2$ is an internal variable, we can remove a contractible 
summand
 $$  0 \lra R\{-1,3\} \xrightarrow{x_4-x_2} R(x_4-x_2)\{-1,1\} \lra 0$$
from the above complex and reduce the ground ring to $R'=\Q[a,x_1,x_4].$ 
We get the complex $(a, x_1-x_4)$ shifted by $\{-1,1\}[-1].$ 
Thus, $C(D_2)\{0,2\}\cong C(\Gamma)\{-1,1\}[-1],$ for $\Gamma$ 
as in figure~\ref{digon}. Since $\Gamma\cong D_1,$ there is an 
isomorphism $C(D_2)\{1,1\}[1]\cong C(D_1).$ We record this as 

\begin{prop} Complexes of matrix factorizations 
 $C(D_1)$ and $C(D_2)\{1,1\}[1]$ are isomorphic as objects of $K(hmf_w),$
with $w=a(x_1-x_4).$  
\end{prop} 

\begin{figure}  \drawing{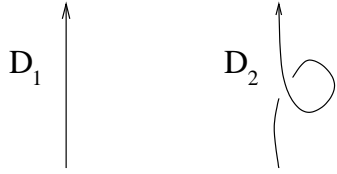}\caption{Type IB Reidemeister move} 
 \label{curl1b} \end{figure}

A similar computation takes care of the Reidemeister move IB: 

\begin{prop} Complexes of matrix factorizations 
  $C(D_1)$ and $C(D_2)$ are isomorphic as objects of $K(hmf_w),$
for $D_1, D_2$ depicted in figure~\ref{curl1b}.
\end{prop} 

\vspace{0.1in}


{\bf 5. Invariance under Reidemeister move IIa.}

 \begin{figure} [htb] \drawing{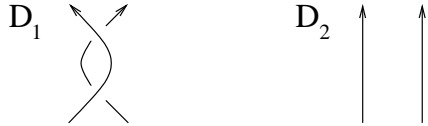}\caption{Type IIa move}
 \label{ta1}
 \end{figure}

 \begin{figure} [htb] \drawing{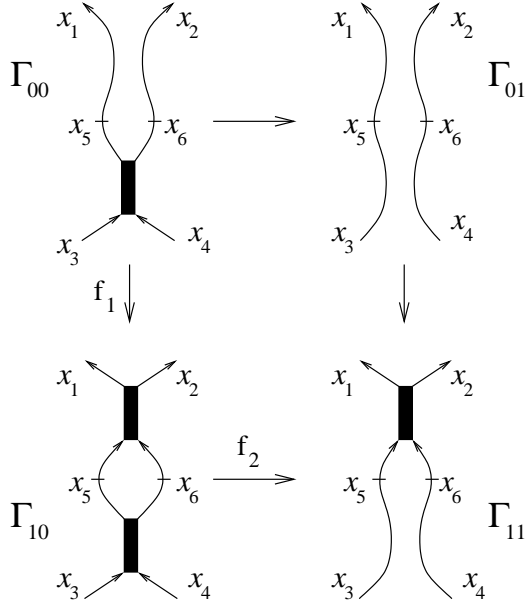}\caption{Four resolutions
 of $D_1$}
 \label{reid2a}
 \end{figure}

Complexes of matrix factorizations $C(D_1)$ and $C(D_2),$ for the diagrams 
depicted in figure~\ref{ta1}, live in the category $K(hmf_w)$ with 
$w=a(x_1+x_2-x_3-x_4),$ viewed as an element of the ground ring 
$R=\Q[a,x_1, x_2,x_3,x_4].$ The complex $C(D_2)\cong C(\Gamma_{01})$ lies 
entirely in cohomological degree zero, since $D_2$ has no crossings. 

\begin{prop} \label{rei-m-t} Complexes of matrix factorizations 
  $C(D_1)$ and $C(D_2)$ are equivalent as objects of $K(hmf_w).$ 
\end{prop} 

\emph{Proof:} It suffices to show that
$f_1$ is an isomorphism (in $hmf_w$) from $C(\Gamma_{00})$ 
to a direct summand of $C(\Gamma_{10}),$ and that there is a decomposition 
$$C(\Gamma_{10}) \cong \mathrm{Im}(f) \oplus M$$ 
with $f_2$ restricting to an isomorphism between $M$ and $C(\Gamma_{11}).$ 
Then $C(D_1)$ would be isomorphic to a direct sum of 
contractible complexes 
\begin{eqnarray*} 
& & 0 \lra C(\Gamma_{00}) \stackrel{\Id}{\lra} \mathrm{Im}(f) \lra 0, \\
& & 0 \lra M \stackrel{f_2}{\lra} C(\Gamma_{11}) \lra 0
\end{eqnarray*} 
and the factorization $C(\Gamma_{01}),$ isomorphic to  $C(D_2).$ 

We start by writing down the diagram of factorizations and maps 
$$ C(\Gamma_{00}) \stackrel{f_1}{\lra} C(\Gamma_{10}) 
\stackrel{f_2}{\lra}C(\Gamma_{11})$$ 
and simplify them in $hmf_w$ by removing contractible direct 
summand factorizations from each of the three terms. The diagram 
has the form 
\begin{eqnarray*} 
& & \left( \begin{array}{cc} a & x_1+ x_2-x_5 -x_6 \\ 0  & x_2-x_6 \\
  a & x_5 + x_6 -x_3 -x_4 \\ 0 & (x_6-x_4)(x_3-x_6) \end{array}
   \right) \stackrel{f_1}{\lra} 
  \left( \begin{array}{cc} a & x_1+ x_2-x_5 -x_6 \\ 0 & (x_2-x_6)(x_5-x_2) \\
  a & x_5 + x_6 -x_3 -x_4 \\ 0 & (x_6-x_4)(x_3-x_6) \end{array}
   \right) \stackrel{f_2}{\lra}  \\
& & \left( \begin{array}{cc} a & x_1+ x_2-x_5 -x_6 \\ 0 & (x_2-x_6)(x_5-x_2) \\
  a & x_5 + x_6 -x_3 -x_4 \\ 0 & x_6-x_4 \end{array}
   \right)
\end{eqnarray*} 
with 
$$f_1= \Id \otimes \psi'(x_5-x_2) \otimes \Id^{\otimes 2}, 
 \hspace{0.2in} 
  f_2= \Id^{\otimes 3} \otimes \psi(x_3-x_6).$$ 
Apply row transformation $[13]_1$ to all three Koszul matrices. 
The new matrices will have identical first rows 
$(a, x_1+x_2-x_3-x_4)$ and identical third rows $(0, x_5+x_6-x_3-x_4).$
We remove the third rows and exclude internal variable $x_5$ substituting 
$x_3+x_4-x_5$ in its place everywhere else. The diagram becomes the 
tensor product of the Koszul factorization $(a, x_1+x_2-x_3-x_4)$  
and the diagram 
\begin{eqnarray*} 
& & 
 \left( \begin{array}{c} x_2-x_6 \\ (x_6-x_4)(x_3-x_6) \end{array} 
  \right) \stackrel{g_1}{\lra} 
  \left( \begin{array}{c} (x_2-x_6)(x_3+x_4-x_6-x_2) \\ (x_6-x_4)(x_3-x_6) 
  \end{array} \right) \stackrel{g_2}{\lra} \\ 
& &  \left( \begin{array}{c} (x_2-x_6)(x_3+x_4-x_6-x_2) \\ x_6-x_4 
  \end{array} \right)
\end{eqnarray*} 
where 
$$ g_1=\psi'(x_3+x_4-x_6-x_2)\otimes \Id, \hspace{0.2in} 
   g_2= \Id \otimes \psi(x_3-x_5).$$ 
We omitted  the first columns from the Koszul matrices, since their terms 
are all zeros. The only internal variable left is $x_6.$ 
The bottom term in the first two factorizations is 
 $$(x_6-x_4)(x_3-x_6) = - x_6^2 + (x_3 + x_4) x_6 -x_3 x_4.$$ 
Let $R'=\Q[a,x_1, x_2, x_3,x_4]$ be the polynomial ring on all 
external variables. Currently we're working over the ring 
$R'[x_6].$ We remove the bottom term from the first two 
factorizations simultaneously reducing to $R',$ imposing 
the relation $x_6^2 = (x_3+x_4)x_6 - x_3 x_4,$ and treating 
multiplication by $x_6$ as an endomorphism of the free $R'$-module 
$R'[x_6]/((x_6-x_4)(x_3-x_6)).$ 
Likewise, in the rightmost factorization, we remove the bottom 
row $(x_6-x_4),$ reduce to the ground ring $R'$ and impose the 
relation $x_6=x_4.$ Our diagram simplifies to 
 \begin{equation*}
   \begin{CD}
     R'1\oplus R'x_6  @>{x_2-x_6}>> R'1\oplus R'x_6   \\
     @V{1}VV      @V{x_3+x_4-x_6-x_2}VV  \\
     R'1\oplus R'x_6 @>{(x_2-x_4)(x_3-x_2)}>> R'1\oplus R' x_6 \\
     @V{x_6\to x_4}VV      @V{x_6\to x_4}VV  \\
     R'    @>{(x_2-x_4)(x_3-x_2)}>> R'  
   \end{CD}
 \end{equation*}
where, for instance, the bottom row denotes the factorization 
$$R' \stackrel{0}{\lra} R' \xrightarrow{ (x_2-x_4)(x_3-x_2)} R' $$
and the maps $g_1, g_2$ are given by vertical arrows.  
Stripping off a contractible summand 
$$ R'1 \stackrel{1}{\lra} R'(x_2-x_6) $$ 
from the first factorization, we reduce it to 
$$ R'(x_6+x_2-x_3-x_4) \xrightarrow{(x_2-x_4)(x_2-x_3)} R' 1. $$ 
The middle factorization is a direct sum of two isomorphic (up to 
grading shift) factorizations 
$$ R' 1  \xrightarrow{(x_2-x_4)(x_3-x_2)} R' 1$$
and 
$$ R'(x_6+x_2-x_3-x_4) \xrightarrow{(x_2-x_4)(x_3-x_2)} 
 R'(x_6+x_2-x_3-x_4).$$ 
The map $g_1$ takes the top factorization (in its reduced form)
isomorphically onto the second summand of the middle factorization. 
The map $g_2$ restricts to an isomorphism from the first direct 
summand of the middle factorization to the bottom factorization. 
Our claim and the proposition follow. $\square$ 

\vspace{0.07in} 

The invariance under the mirror image of the figure~\ref{ta1} move 
can be verified similarly. 

\vspace{0.1in}


{\bf 6. Invariance under Reidemeister move III.} 

\vspace{0.1in} 

Let factorization $\Upsilon$ be given by the following Koszul matrix: 
\begin{equation*} 
 \left( \begin{array}{cc} a & x_1+x_2+x_3-x_4-x_5-x_6 \\
       0  & x_1x_2+x_1x_3+x_2x_3 - x_4x_5 -x_4x_6 -x_5x_6 \\
       0  & x_1 x_2 x_3 -x_4 x_5 x_6 \end{array} \right) 
\end{equation*} 
The gradings are normalized so that the differential has bidegree $(1,1).$ 
For instance, the bottom row denotes the factorization 
$$R \stackrel{0}{\lra} R\{-1,5\} \xrightarrow{x_1 x_2 x_3 -x_4 x_5 x_6}
       R $$
with $R=\Q[a,x_1, \dots, x_6].$ The potential is $w=a(x_1+x_2+x_3-x_4-x_5-x_6).$
                
\begin{figure}  \drawing{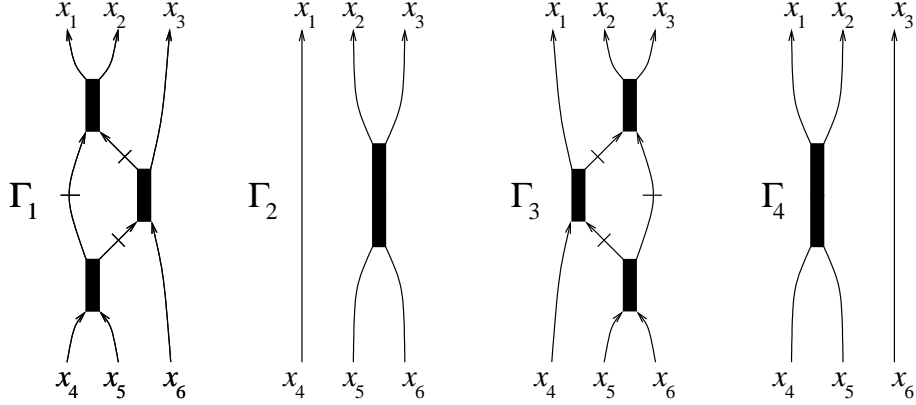}\caption{Diagrams $\Gamma_1, \Gamma_2, 
\Gamma_3, \Gamma_4$} 
 \label{fourgr} \end{figure}

\begin{prop}\label{t-three} In $hmf_w$ there are isomorphisms 
\begin{eqnarray} 
   C(\Gamma_1) & \cong & C(\Gamma_4)\{0,2\} \oplus  \Upsilon ,  \\
 C(\Gamma_3) & \cong & C(\Gamma_2)\{0,2\} \oplus  \Upsilon , 
\end{eqnarray} 
for $\Gamma_1, \Gamma_2, \Gamma_3, \Gamma_4$ depicted in figure~\ref{fourgr}.  
\end{prop} 

\emph{Proof:} To prove the first isomorphism, we place labels $x_7, x_8, x_9$ 
(from top to bottom) on the three marks of $\Gamma_1$ and write 
$C(\Gamma_1)$ in the Koszul form 
$$ C(\Gamma_1) = \left( \begin{array}{cc} 
    a & x_ 1+ x_2 - x_8  - x_7 \\ 
    0 & (x_2-x_7)(x_8-x_2)  \\
    a & x_7 + x_3 - x_9 - x_6 \\
    0 & (x_3-x_6)(x_9-x_3) \\
    a & x_8 + x_9 -x_4 -x_5 \\
    0 & (x_9-x_5)(x_4 - x_9) \end{array} \right) $$ 
Applying transformations $[13]_1$ and $[15]_1$ we get a matrix 
with the third and fifth rows 
$$ (0, x_7 + x_3 -x_9 -x_6), \hspace{0.2in} (0, x_8+x_9 -x_4 -x_5).$$ 
We use these rows to exclude internal variables $x_7$ and $x_8$ and 
reduce $C(\Gamma_1)$ to the following Koszul form 
$$ \left( \begin{array}{cc} 
    a & x_ 1+ x_2 +x_3 - x_4  - x_5 -x_6 \\ 
    0 & (x_2+x_3-x_6-x_9)(x_4+x_5-x_2-x_9)  \\
    0 & (x_3-x_6)(x_9-x_3) \\
    0 & (x_9-x_5)(x_4 - x_9) \end{array} \right) $$ 
Notice that variables $a$ and $x_1$ appear only in the first row. 
Moreover, Koszul forms of factorizations $\Upsilon$ and 
$C(\Gamma_4)$ have the same first row. Next, we ignore 
the first row of $C(\Gamma_1)$ and operate on the other three 
rows. The first column of the Koszul matrix then consists of zeros 
and we omit it. To simplify the factorization
$$ \left( \begin{array}{c}
    (x_2+x_3-x_6-x_9)(x_4+x_5-x_2-x_9)  \\
    (x_3-x_6)(x_9-x_3) \\
    (x_9-x_5)(x_4 - x_9) \end{array} \right) $$ 
we use the last term to reduce to at most linear terms in the 
last remaining internal variable $x_9.$ Remove the last row 
and impose the relation $x_9^2=(x_4+x_5)x_9 -x_4 x_5.$ 
Modulo this relation and after adding the second row, the 
first row loses $x_9$ and the matrix becomes
$$\left( \begin{array}{c} 
  (x_3-x_6)(x_4+x_5-x_2-x_3)+(x_2-x_5)(x_4-x_2) \\
    (x_3-x_6)(x_9-x_3) \end{array}\right) $$ 
Now $x_9$ appears only in the bottom row, which we can 
write as 
$$R 1 \oplus R x_9 \xrightarrow{ (x_3-x_6)(x_9-x_3)} 
  R 1 \oplus R x_9.$$
Changing basis of the free $R$-module on the left hand side 
from $\{1, x_9\}$ to $\{1, x_9+x_3-x_4-x_5\}$ 
and of the module on the right to $\{ 1, x_9-x_3\},$ we decompose 
this complex into a direct sum of 
$$ R (x_9+x_3-x_4-x_5) \xrightarrow{(x_3-x_4)(x_5-x_3)(x_3-x_6)} 
 R1 $$ 
and 
$$ R 1 \xrightarrow{x_3-x_6} R(x_9-x_3).$$ 
Adding the other rows, we obtain a decomposition of $C(\Gamma_1)$ 
into direct sum of factorizations with Koszul matrices 
$$ \left( \begin{array}{cc} 
    a & x_ 1+ x_2 +x_3 - x_4  - x_5 -x_6 \\ 
    0 & (x_3-x_6)(x_4+x_5-x_2-x_3)+(x_2-x_5)(x_4-x_2)  \\
    0 & (x_3-x_4)(x_3-x_5)(x_3-x_6) \end{array} \right) $$ 
and 
$$\left( \begin{array}{cc} 
    a & x_ 1+ x_2 +x_3 - x_4  - x_5 -x_6 \\ 
    0 & (x_3-x_6)(x_4+x_5-x_2-x_3)+(x_2-x_5)(x_4-x_2)  \\
    0 & x_3-x_6 \end{array} \right), $$
the latter shifted by $\{0,2\}$ due to the bidegree 
$(0,2)$ vector $x_9-x_3$ being a generator of the module $R(x_9-x_3).$ 
It is easy to check that the matrices above describe factorizations 
$\Upsilon$ and $C(\Gamma_4),$ respectively. $\square$ 

\vspace{0.1in} 

\begin{figure}  \drawing{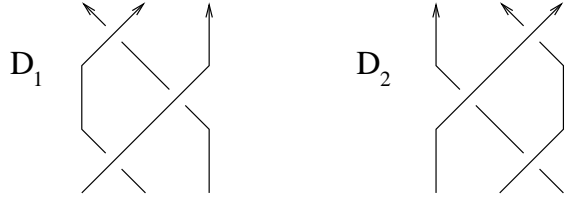}\caption{Reidemeister move III} 
 \label{rthree0} \end{figure}

\begin{prop} \label{inv-r-three} 
Complexes $C(D_1)$ and $C(D_2),$ for diagrams depicted  
in figure~\ref{rthree0}, are isomorphic in the category $K(hmf_w).$ 
\end{prop}  

\begin{figure}  \drawing{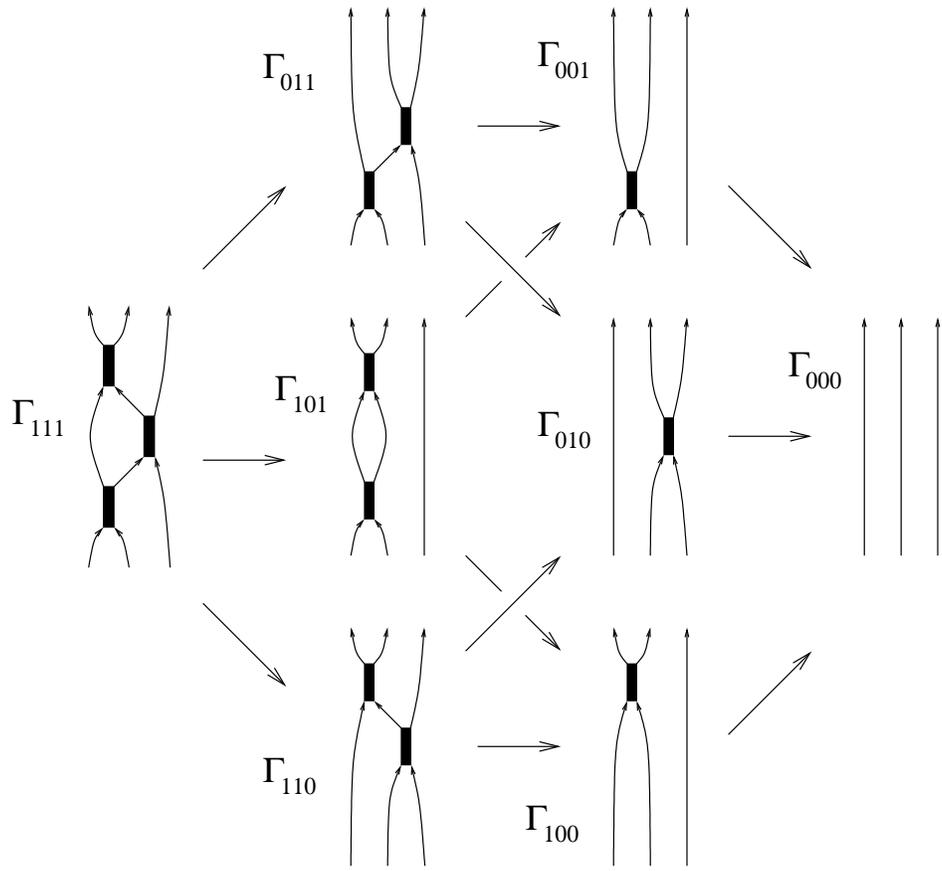}\caption{Resolution cube of $D_1$} 
 \label{rthree1} \end{figure}

\emph{Proof} is similar to the one in [KR]. The complex $C(D_1)$ 
consists of eight factorizations assigned to  diagrams depicted in 
figure~\ref{rthree1} (also see figure~\ref{cones}). We ignore the 
overall shift by $\{0,-2\}$ in the resolution of each crossing, which was 
needed for the invariance under the Reidemeister move IIa, but does not
make any difference for Reidemeister move III. Proposition~\ref{t-three} 
tells tells us that 
$$C(\Gamma_{111}) \cong \Upsilon \oplus C(\Gamma_{100})\{0,2\}$$ 
(also observe that $\Gamma_{100}\cong \Gamma_{001}$),  
while our proof of the invariance under the Reidemeister move IIa implies 
$$ C(\Gamma_{101})\cong C(\Gamma_{100})\{0,2\} \oplus  C(\Gamma_{100}).$$ 
A computation similar to the one in that proof shows that the map 
$\chi_1: C(\Gamma_{111}) \lra C(\Gamma_{101}),$ when 
restricted to the direct summand  isomorphic to $C(\Gamma_{100})\{0,2\},$
is an isomorphism onto a direct summand of $C(\Gamma_{101}),$ while  
our proof of proposition~\ref{rei-m-t} implies that 
$\chi_1:C(\Gamma_{101}) \lra C(\Gamma_{001})$ is an isomorphism
when restricted to the direct summand $C(\Gamma_{100})$ of 
$C(\Gamma_{101}).$ 

\begin{figure}  \drawing{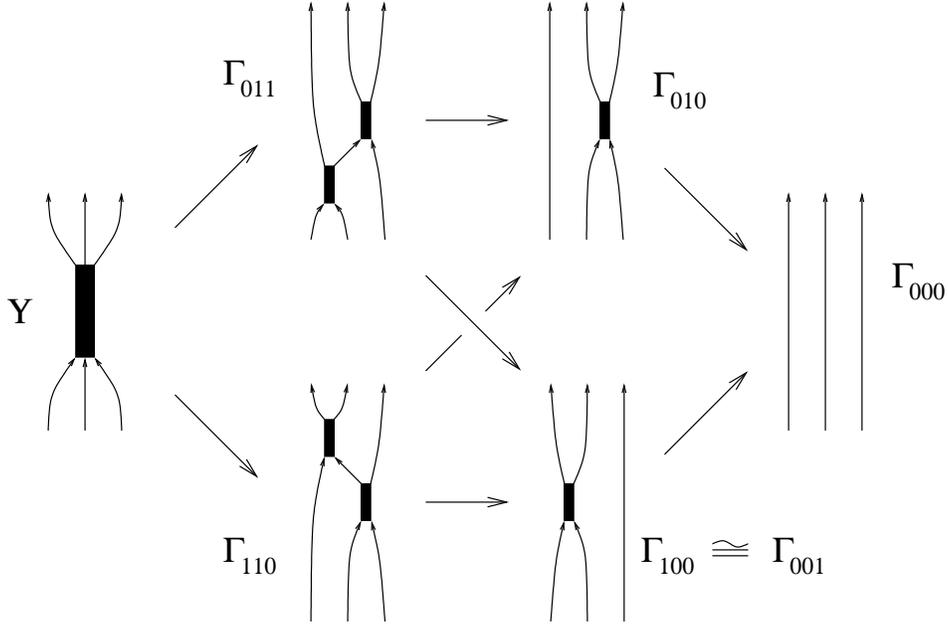}\caption{Complex $C'$} 
 \label{rthree2} \end{figure}

After removing contractible summands 
$$ 0 \lra C(\Gamma_{100}) \{ 0,2\} \stackrel{\cong}{\lra} 
       C(\Gamma_{100}) \{ 0,2\} \lra 0 $$ 
and
$$0 \lra C(\Gamma_{100})  \stackrel{\cong}{\lra} 
       C(\Gamma_{100})  \lra 0$$
we reduce $C(D_1)$ to a complex $C'$ that is graphically 
depicted in figure~\ref{rthree2}, factorization $\Upsilon$ 
assigned to the diagram $Y.$

\begin{lemma} \label{lemma-indec} 
The complex $C'\cong C(D_1)$ is indecomposable in the 
category $K(hmf_w).$ 
\end{lemma} 

In other words, we cannot write $C'\cong M\oplus N$ for two nontrivial 
objects $M,N$ of $K(hmf_w).$ Indeed, invariance under the Reidemeister 
move IIa tells us that tensoring with a complex of factorizations assigned 
to a crossing is an invertible functor. Precisely, it's an invertible 
functor from the category $K(hmf_v)$ to $K(hmf_u)$ where 
$v=ax_3+ax_4+f(\mathbf{x}),$ 
 $u=ax_1+ax_2 + f(\mathbf{x}),$ and $f(\mathbf{x})$ is  any polynomial 
in variables $\mathbf{x}$ disjoint from $x_1, \dots, x_4.$ An invertible 
functor is indecomposable iff the identity functor is. 
The identity functor, in general, corresponds 
to the diagram comprised of $n$ parallel lines, compatibly
oriented (the diagram of the trivial braid). Its factorization $S$  
can be written as the tensor product of 
$(a, x_i -x_{n+i}), $ over $i=1, \dots, n.$ An easy computation 
(for instance, as in the proof of the next lemma) 
shows that the hom space $\Hom_{hmf_w}(S,S)$ of bigrading-preserving 
factorization homomorphisms up to chain homotopies is one-dimensional. 
Therefore, $S$ and the identity functor are indecomposable, for 
otherwise a projection onto a direct summand would ensure that the above 
hom space is at least 2-dimensional. $\square$

\begin{lemma} \label{space-one-dim} 
For any arrow $\Gamma\lra \Gamma'$ in figure~\ref{rthree2}, the space 
of bidegree-preserving maps $C(\Gamma) \to C(\Gamma')$ is 
one-dimensional (over the ground field $\Q$) and is generated by $\chi_1.$ 
\end{lemma}

We can prove the lemma on a case-by-case basis, separately for each arrow. 
In general, to compute the dimension of $\Hom_{hmf_w}(M,N),$ for matrix 
factorizations $M,N,$ with $N$ of finite rank, we use the isomorphism 
$$ \mathrm{EXT}_{hmf_w}(M,N) \cong H(N \otimes_R M_{\bullet})$$ 
where EXT refers to taking ext groups of the pair $M,N$ in all bidegrees, 
$M_{\bullet}$ is the $R$-module dual of $M,$ and $H$ stands for cohomology. 
Restricting the left hand side to $\Hom$ corresponds to taking the bidegree 
$(0,0)$ summand of the right hand side. The dual of a Koszul factorization 
$(a_i, b_i)$ is the Koszul factorization $(b_i, -a_i),$ with suitably 
shifted gradings. 

For instance, to determine the dimension of the space 
$$\Hom_{hmf_w}(C(\Gamma_{110}), C(\Gamma_{100}))$$ (the bottom 
arrow in figure~\ref{rthree2}) 
we first write the Koszul matrix of $C(\Gamma_{110}):$
$$ \left( \begin{array}{cccc}
 \{0,0\} & a &  \{1,-1\} & x_1+x_2-x_4-x_7  \\
 \{0,0\} & a &  \{1,-1\} & x_7+x_3-x_5-x_6  \\
 \{0,0\} & 0 &  \{1, -3\} & (x_2-x_4)(x_2-x_7)  \\
 \{0,0\} & 0 &  \{1,-3\} & (x_3-x_5)(x_3-x_6) \end{array} \right) $$ 
Here $x_7$ is the variable assigned to the internal mark of 
$\Gamma_{110}.$  We also added two columns indicating the bidegrees of $R.$ 
For instance, the second row denotes the factorization
$$ R\{0,0\} \stackrel{a}{\lra} R\{1,-1\} 
  \xrightarrow{x_7+x_3-x_5-x_6} R\{0,0\}.$$
After we apply $[12]_1,$ the second row becomes $(0,x_7+x_3-x_5-x_6);$ 
 we get rid of it and of the variable $x_7.$ Thus, $C(\Gamma_{110})$ 
is isomorphic to the factorization assigned to the Koszul matrix 
$$  \left( \begin{array}{cccc}
 \{0,0\} & a &  \{-1,1\} & x_1+x_2+x_3-x_4-x_5-x_6  \\
 \{0,0\} & 0 &  \{-1, 3\} & (x_2-x_4)(x_2+x_3-x_5-x_6)  \\
 \{0,0\} & 0 &  \{-1,3\} & (x_3-x_5)(x_3-x_6) \end{array} \right). $$ 
The dual $C(\Gamma_{110})_{\bullet}$ of $C(\Gamma_{110})$ can 
be represented by the matrix 
$$  \left( \begin{array}{cccc}
 \{0,0\} & x_1+x_2+x_3-x_4-x_5-x_6 &  \{1,-1\} & -a  \\
 \{0,0\} & (x_2-x_4)(x_2+x_3-x_5-x_6) &  \{1,-3\} & 0  \\
 \{0,0\} &(x_3-x_5)(x_3-x_6) &  \{1,-3\} & 0 \end{array} \right) $$ 
and the tensor product complex $C(\Gamma_{110})_{\bullet}\otimes 
C(\Gamma_{100})$ by 
$$ \left( \begin{array}{cccc}
\{0,0\} & a & \{-1,1\} & x_1+x_2-x_4-x_5  \\
 \{0,0\} & 0 & \{-1,3\} &  (x_2-x_4)(x_2-x_5) \\
  \{0,0\} & a & \{-1,1\} &  x_3-x_6 \\
 \{0,0\} & x_1+x_2+x_3-x_4-x_5-x_6 &  \{1,-1\} & -a  \\
 \{0,0\} & (x_2-x_4)(x_2+x_3-x_5-x_6) &  \{1,-3\} & 0  \\
 \{0,0\} &(x_3-x_5)(x_3-x_6) &  \{1,-3\} & 0 \end{array} \right), $$
where the first three rows describe $C(\Gamma_{100}).$  
We do transformation $[13]_1$ and shift rows 4 and 5 by one each. 
We get 
$$ \left( \begin{array}{cccc}
 \{0,0\} & a & \{-1,1\} & x_1+x_2+x_3-x_4-x_5-x_6  \\
 \{0,0\} & 0 & \{-1,3\} & (x_2-x_4)(x_2-x_5) \\
 \{0,0\} & 0 & \{-1,1\} & x_3-x_6 \\
 \{1,-1\} & -a & \{0,0\} &   x_1+x_2+x_3-x_4-x_5-x_6 \\
 \{1,-3\} & 0 & \{0,0\} & (x_2-x_4)(x_2+x_3-x_5-x_6)   \\
 \{0,0\} &(x_3-x_5)(x_3-x_6) &  \{1,-3\} & 0 \end{array} \right) $$ 
We apply the transformation $[14]_{-1},$ then shift rows $1$ and $6$ 
to obtain  
$$ \left( \begin{array}{cccc}
 \{-1,1\} & 0 & \{0,0\} & a  \\
 \{0,0\} & 0 & \{-1,3\} & (x_2-x_4)(x_2-x_5) \\
 \{0,0\} & 0 & \{-1,1\} & x_3-x_6 \\
 \{1,-1\} & 0 & \{0,0\} &   x_1+x_2+x_3-x_4-x_5-x_6 \\
 \{1,-3\} & 0 & \{0,0\} & (x_2-x_4)(x_2+x_3-x_5-x_6)   \\
 \{1,-3\} & 0 & \{0,0\} & (x_3-x_5)(x_3-x_6) \end{array} \right). $$ 
Cross out row $3$ and convert $x_6$ to $x_3$ everywhere else. The 
matrix reduces to 
$$ \left( \begin{array}{cccc}
 \{-1,1\} & 0 & \{0,0\} & a  \\
 \{0,0\} & 0 & \{-1,3\} & (x_2-x_4)(x_2-x_5) \\
 \{1,-1\} & 0 & \{0,0\} &   x_1+x_2-x_4-x_5 \\
 \{1,-3\} & 0 & \{0,0\} & (x_2-x_4)(x_2-x_5)   \\
 \{1,-3\} & 0 & \{0,0\} & 0 \end{array} \right). $$ 
We apply $[24]_1,$ then remove rows $1$ and $3$ simultaneously 
with getting rid of the variables $a$ and $x_5.$ The resulting matrix is 
$$\left( \begin{array}{cccc}
 \{0,0\} & 0 & \{-1,3\} & (x_2-x_4)(x_4-x_1) \\
 \{1,-3\} & 0 & \{0,0\} & 0  \\
 \{1,-3\} & 0 & \{0,0\} & 0 \end{array} \right). $$
Let now $R'=\Q[x_1,x_2,x_3,x_4].$ The cohomology of the complex 
described by this matrix is the tensor product of 
the quotient $R'/((x_2-x_4)(x_4-x_1))$ and the bigraded 
vector space
$$(\Q\{1,-3\} \oplus \Q )\otimes (\Q\{1,-3\} \oplus \Q ).$$ 
The bigraded dimension of $R'/((x_2-x_4)(x_4-x_1))$ has 
the form $1+\alpha$ where $\alpha\in q^2 \Z[q^2],$ while 
that of the second term is 
$(1+ tq^{-3})^2.$ Therefore, the bigraded dimension of 
the cohomology of the complex $C(\Gamma_{110})_{\bullet}\otimes 
C(\Gamma_{100})$ has the form 
$$(1+\alpha)(1+ 2t q^{-3}+ t^2 q^{-6}).$$
Writing it as a polynomial in $t$ with coefficients being power 
series in $q,$ we see that the coefficient of the term $t^0 q^0$ 
equals $1.$ Therefore, the bidegree $(0,0)$ summand of the homology 
is one-dimensional, and the hom space 
$\Hom_{hmf_w}(C(\Gamma_{110}), C(\Gamma_{100}))$
has dimension $1.$ 

To show that $\chi_1: C(\Gamma_{110})\lra 
C(\Gamma_{100})$ (corresponding to the splitting of the right wide 
edge of $\Gamma_{110}$ into two parallel lines) 
 generates this one-dimensional space, it suffices to show that 
$\chi_1$ is not null-homotopic. We can write the factorizations and 
the map in the following Koszul form 
$$ \left( \begin{array}{cc} a & x_1+x_2-x_4-x_7  \\ 
  0 & (x_2-x_4)(x_2-x_7)  \\
  a & x_7+x_3-x_5-x_6  \\
  0 & (x_3-x_5)(x_3-x_6) \end{array} \right) \stackrel{\chi_1}{\lra} 
  \left( \begin{array}{cc} a & x_1+x_2-x_4-x_7  \\ 
  0 & (x_2-x_4)(x_2-x_7)  \\
  a & x_7+x_3-x_5-x_6  \\
  0 & x_3-x_6 \end{array} \right) $$ 
with $\chi_1=\Id^{\otimes 3}\otimes \psi(x_3-x_5).$ Applying 
row transformation $[13]_1$ to each matrix and then excluding $x_7$ 
we reduce the map to 
$$ \left( \begin{array}{cc} a & x_1+x_2+x_3-x_4-x_5-x_6  \\ 
  0 & (x_2-x_4)(x_2+x_3-x_5-x_6)  \\
  0 & (x_3-x_5)(x_3-x_6) \end{array} \right) \stackrel{\chi_1}{\lra} 
  \left( \begin{array}{cc} a & x_1+x_2+x_3-x_4-x_5-x_6  \\ 
  0 & (x_2-x_4)(x_2+x_3-x_5-x_6)  \\
  0 & x_3-x_6 \end{array} \right) $$ 
with $\chi_1=\Id^{\otimes 2}\otimes \psi(x_3-x_5).$ 
Turn both $\Gamma_{110}$ and $\Gamma_{100}$ into closed 
diagrams $\widehat{\Gamma}_{110}$ and $\widehat{\Gamma}_{100}$
by connecting top endpoints of each diagram with its bottom endpoints
by three disjoint arcs. To check that $\chi_1$ is not null-homotopic, 
it's enough to verify that the induced map on cohomology 
$$\widehat{\chi}_1: H(\widehat{\Gamma}_{110}) \lra 
  H(\widehat{\Gamma}_{100})$$
is non-trivial. We represent this map in Koszul form as 
$$ \left( \begin{array}{cc} a & x_1+x_2+x_3-x_4-x_5-x_6  \\ 
  0 & (x_2-x_4)(x_2+x_3-x_5-x_6)  \\
  0 & (x_3-x_5)(x_3-x_6) \\
  a & x_4 -x_1 \\
  a & x_5 -x_2 \\
  a & x_6 -x_3 \end{array} \right) \stackrel{\widehat{\chi}_1}{\lra} 
  \left( \begin{array}{cc} a & x_1+x_2+x_3-x_4-x_5-x_6  \\ 
  0 & (x_2-x_4)(x_2+x_3-x_5-x_6)  \\
  0 & x_3-x_6   \\
  a & x_4 -x_1 \\
  a & x_5 -x_2 \\
  a & x_6 -x_3
 \end{array} \right), $$ 
and $\widehat{\chi}_1= \Id^{\otimes 2}\otimes \psi(x_3-x_5) \otimes 
\Id^{\otimes 3}.$ Doing transformations $[14]_1,$ $[15]_1,$ $[16]_1$ 
and excluding $x_4,$ $x_5,$ $x_6$ and $a,$ we reduce the map 
to the form 
$$ \left( \begin{array}{cc} 0 & 0 \\
  0 & 0 
 \end{array} \right) \stackrel{\widehat{\chi}_1}{\lra} 
  \left( \begin{array}{cc} 0 & 0 \\
  0 & 0
 \end{array} \right), $$ 
where $\widehat{\chi}_1= \Id \otimes \psi(x_3-x_2)$ and the ground 
ring is $\Q[x_1,x_2,x_3]$ (we took the quotient by the 
ideal $(x_1-x_4, x_2-x_5, x_3-x_6,a)$). Clearly, 
$\widehat{\chi}_1$ induces a nontrivial map on 
cohomology, and $\chi_1$ is not null-homotopic.  

\vspace{0.06in} 

Using symmetries of the graphs and factorizations, the other 
cases of the lemma can be reduced to verifying that the hom spaces 
$$\Hom_{hmf_w}(\Upsilon, C(\Gamma_{110}))\hspace{0.2in} \mathrm{and}
\hspace{0.2in}\Hom_{hmf_w}(C(\Gamma_{100}), C(\Gamma_{000}))$$
are both one-dimensional and generated by $\chi_1.$ 
Actual computations, similar to the 
one above, are left to a curious reader. For the first of the 
two hom spaces, by $\chi_1$ we mean the composition of 
$\chi_1: C(\Gamma_{111}) \lra C(\Gamma_{110})$ 
with the inclusion of $\Upsilon$ as a direct 
summand of $C(\Gamma_{111}).$ 
$\square$ 

\vspace{0.1in} 

\begin{lemma} \label{not-iso} For any arrow $\Gamma\lra \Gamma'$ 
in figure~\ref{rthree2}, factorizations $C(\Gamma)$ and $C(\Gamma')$ 
are not isomorphic in $hmf_w.$ \end{lemma} 

\emph{Sketch of proof:} Form the closures $\widehat{\Gamma}$ and $\widehat{\Gamma}'$ 
by connecting top endpoints of each diagram with its bottom endpoints 
by 3 disjoint arcs. A direct computation shows that complexes 
$C(\widehat{\Gamma})$ and $C(\widehat{\Gamma}')$ have non-isomorphic 
cohomology groups (their two-variable Poincare polynomials are different). 
$\square$ 

\vspace{0.1in} 

Thus, the complex $C',$ depicted in figure~\ref{rthree2},
consists of 6 factorizations and its differential is a sum of 
10 maps, one for each arrow of the figure. Each map is either 0 
or a nonzero multiple of the unique (up to rescaling) nontrivial 
map between the two factorizations. For each arrow $b: \Gamma\lra \Gamma'$ 
choose a nontrivial map $m(b): C(\Gamma) \to C(\Gamma').$ 

\begin{lemma} \label{2-nontriv} 
For any two composable arrows $\Gamma\stackrel{b_1}{\lra} \Gamma'
\stackrel{b_2}{\lra}\Gamma''$ the composition $m(b_2)m(b_1)$ is nontrivial 
in $hmf_w.$ 
\end{lemma}

\emph{Proof:} It suffices to check that the composition 
\begin{equation} \label{eq-4-maps} 
\Upsilon \subset C(\Gamma_{111}) \stackrel{\chi_1}{\lra} C(\Gamma_{110}) 
\stackrel{\chi_1}{\lra} C(\Gamma_{100}) \stackrel{\chi_1}{\lra} C(\Gamma_{000})
\end{equation} 
is not null-homotopic. Denote the composition of the last 3 maps by $\chi'_1$ and 
the corresponding "adjoint" composition 
$$ C(\Gamma_{000})\stackrel{\chi_0}{\lra}C(\Gamma_{100})
\stackrel{\chi_0}{\lra}C(\Gamma_{110}) 
\stackrel{\chi_0}{\lra}C(\Gamma_{111})$$
by $\chi'_0.$ We claim that the map $\chi'_1\mathrm{pr} \chi'_0$ is non-zero, 
where $\mathrm{pr}$ is the projection from $C(\Gamma_{111})$ onto 
its direct summand $\Upsilon.$ The map $\chi'_0$ has degree $(0,6)$
and the product $\chi'_1\chi'_0$ is equal to the multiplication by 
$(x_4-x_2)^2(x_5-x_3)$ endomorphism of $C(\Gamma_{000}),$ since 
the composition $\chi_1\chi_0$ is the multiplication by a suitable linear 
combination of $x$'s. The complementary direct summand of 
$C(\Gamma_{111})$ is isomorphic to $C(\Gamma_{100})\{0,2\}.$ 
Denote by $\widetilde{\mathrm{pr}}$ the projection onto this direct summand. Then 
$\mathrm{pr} + \widetilde{\mathrm{pr}}$ is the identity endomorphism of 
$C(\Gamma_{111}).$ 

The composition $\chi'_1\widetilde{\mathrm{pr}} \chi'_0$ factors though a 
degree $(0,2)$ endomorphism of $C(\Gamma_{110}).$ This endomorphism 
is a composition 
$$ C(\Gamma_{110}) \lra C(\Gamma_{100}) \lra C(\Gamma_{110})$$ 
where the first map has degree $(0,0)$ and the second--degree $(0,2).$ These 
maps are, necessarily, rational multiples of $\chi_0$ and $\chi_1$ (corresponding 
to the right wide edge of $\Gamma_{110}$) and their composition  is 
a rational multiple of the multiplication by $x_3-x_5.$ Hence, the composition 
$\chi'_1\widetilde{\mathrm{pr}} \chi'_0$ is a rational multiple of the multiplication 
by $(x_3-x_5)^2 (x_2-x_4).$ To show that 
$$\chi'_1\mathrm{pr} \chi'_0=\chi'_1\chi'_0 - \chi'_1\widetilde{\mathrm{pr}}\chi'_0$$
is not null-homotopic, we observe that the right hand side is the multiplication by 
$$(x_4-x_2)^2(x_5-x_3) - \mu (x_3-x_5)^2 (x_2-x_4)$$
endomorphism of $C(\Gamma_{000}),$ for some rational $\mu.$ The image of 
$\Q[x_1, \dots, x_6]$ in the 
endomorphism ring of $C(\Gamma_{000})$ is the quotient ring by relations 
$x_1=x_4,$ $x_2=x_5$ and $x_3=x_6.$ The polynomial above simplifies to 
$$ (x_1-x_2)^2 (x_2-x_3) - \mu (x_3-x_2)^2 (x_2-x_1) \not= 0 $$ 
in $\Q[x_1, x_2, x_3].$ Therefore, the composition   
$\chi'_1\mathrm{pr} \chi'_0$ is not null-homotopic, and so is the map 
$\Upsilon\lra C(\Gamma_{000})$ in formula (\ref{eq-4-maps}). 
Lemma~\ref{2-nontriv} follows. $\square$ 

\vspace{0.15in} 

The differential in the complex $C'$ can be written as 
$$ d= \sum_b  \lambda_b  m(b), $$
with $\lambda_b\in \Q,$ and  the sum over all arrows $b.$   

\begin{lemma}\label{all-nonzero} All coefficients $\lambda_b$ 
are nonzero rational numbers.  \end{lemma} 

Assume otherwise: $\lambda_b=0$ for some $b.$ Every square in 
the diagram of $C'$ anticommutes, and from lemma~\ref{2-nontriv} 
we derive that some other $\lambda$'s would have to be zero. 
In fact, there will be enough zero maps to split the complex 
into the direct sum of at least two subcomplexes, each comprised 
of two or four factorizations in figure~\ref{rthree2}.
Specifically, the complex will either decompose into a direct sum 
of 3 subcomplexes  of the form 
\begin{equation} \label{sub-com-2}  
0 \lra C(\Gamma) \stackrel{m(b)}{\lra} C(\Gamma') \lra 0
\end{equation}  
for some three arrows $b,$ or as the direct sum of one subcomplex 
of type (\ref{sub-com-2}) and the complementary summand containing 
the other four factorizations. 

A decomposition of $C'$ into a direct sum of 3 subcomplexes 
contradicts lemmas~\ref{lemma-indec}, \ref{not-iso}. To see the 
impossibility of the decomposition of 
the second kind, it's enough to show that the complementary summand 
cannot be trivial in $hmf_w.$ This summand would consist of four factorizations 
that sit in the vertices of one of the four squares in figure~\ref{rthree2}. 
For instance, it could have the form 
$$ 0 \lra C(Y) \lra C(\Gamma_{110})\oplus C(\Gamma_{011}) \lra 
  C(\Gamma_{010}) \lra 0.$$ 
Triviality of the summand would imply that its identity map is null-homotopic. 
In particular, the identity map of the rightmost factorization in the complex would 
factor through a map to the middle term. This map should have bidegree $(0,0).$ 
The following lemma establishes the contradiction. 

\begin{lemma} For any arrow $\Gamma \lra \Gamma'$ in figure~\ref{rthree2} 
we have  
$$\Hom_{hmf_w}(C(\Gamma'), C(\Gamma)) = 0.$$
\end{lemma} 
Thus, any bidegree zero map is trivial. The lemma can be proved in the 
same way as lemma~\ref{space-one-dim}. $\square$ 

\vspace{0.1in} 

Lemma~\ref{all-nonzero} follows. $\square$ 

\vspace{0.1in} 

To summarize, we established that the coefficients $\lambda_b$ 
in the differential for the complex $C'$ are all nonzero. Rescaling, if necessary, 
we can turn them into $1$'s and $-1$'s. Moreover, the complex $C'$ is uniquely 
determined, up to isomorphism, by the condition that $\lambda_b\not= 0$ for 
all $b.$ We have $C(D_1)\cong C'.$ Nearly identical arguments show that $C(D_2) \cong C'$ 
as well. Therefore, $C(D_1) \cong C(D_2),$ and proposition~\ref{inv-r-three} follows. $\square$ 

\vspace{0.15in}


{\bf 7. Computing the Euler characteristic.} 

\vspace{0.1in} 

Given a braid diagram $D,$ the complexes $C(D\sigma^{\pm 1}_i)$ are, 
up to shifts, the cones of maps $\chi_0, \chi_1$ between factorizations 
assigned to diagrams $D$ and $De_i$ where $e_i$ denotes a wide edge 
placed between $i$-th and $(i+1)$-st strands of the braid. See 
Section~\ref{section-one} and figure~\ref{cones} for details. Recall that 
$\langle D\rangle$ denotes the Euler characteristic of cohomology $H(D).$ 
Our definition of $C(D\sigma^{\pm 1}_i)$ via cones implies the following 
relations on the Euler characteristics: 
\begin{eqnarray*} 
 \langle D\sigma_i \rangle & = & \langle De_i \rangle - q^2 \langle  D \rangle , \\
 \langle D\sigma^{-1}_i \rangle & = & q^{-2}( \langle De_i \rangle - 
   \langle  D \rangle).
 \end{eqnarray*} 
Excluding $\langle De_i \rangle,$ we get 
\begin{equation*} 
  q^{-1} \langle D\sigma_i \rangle - q^{-1}  \langle D\sigma^{-1}_i \rangle 
    = ( q^{-1} - q)  \langle  D \rangle ,
\end{equation*} 
which is the relation satisfied by the function $F(D)$ defined in Section~\ref{section-one}. 
If $D$ is the one-strand braid diagram of the unknot, $H(D) \cong \Q[x]\{-1,1\}$ 
and 
$$\langle D \rangle = \frac{t^{-1}}{q^{-1}-q},$$ 
which is our normalization of $F(D).$ Propositions proved above imply that 
$\langle D \rangle$ satisfies all other defining properties of $F(D).$ Therefore, 
$\langle D \rangle = F(D)$ and Theorem~\ref{sec-thm} follows. $\square$

\vspace{0.2in} 

{\bf Acknowledgements:} The idea to consider matrix factorizations with 
the potential $ax$ came up when the authors were attending Virginia 
Topology Conference in December of 2004. We'd like to thank 
Slava Krushkal and Frank Quinn, the conference organizers, 
for bringing 
us together and creating an inspiring research atmosphere. 
M.K. is grateful to Sergei Gukov, 
Peter Ozsv\'ath, Jacob Rasmussen, Albert Schwarz and Oleg Viro for 
interesting discussions. We'd like to acknowledge 
NSF support via grants DMS-0407784 
and DMS-0196131.

\vspace{0.2in} 

{\bf References} 

\noindent
[APS] M.~Asaeda, J.~Przytycki and A.~Sikora, Categorification of the 
Kauffman bracket skein module of I-bundles over surfaces, Algebr. Geom. Topol. 
4 (2004) 1177-1210, arXiv  math.QA/0409414. 

\noindent 
[AF] B.~Audoux and T.~Fiedler, A Jones polynomial for braid-like 
isotopies of oriented links and its categorification, arXiv math.GT/0503080. 

\noindent 
[BN] D.~Bar-Natan, Khovanov's homology for tangles and cobordisms, 
arXiv math.GT/0410495. 

\noindent 
[G] B.~Gornik, Note on Khovanov link cohomology, arXiv 
 math.QA/0402266.  

\noindent
[GSV] S.~Gukov, A.~Schwarz and C.~Vafa, Khovanov-Rozansky homology 
and topological strings, arXiv hep-th/0412243. 

\noindent 
[HOMFLY] P.~Freyd, D.~Yetter, J.~Hoste, W.~B.~R.~Lickorish, 
K.~Millett and A.~Ocneanu, A new polynomial invariant of knots 
and links, \emph{Bull. AMS. (N.S.)} {\bf 12} 2, 239--246, 1985. 

\noindent 
[K] M.~Khovanov, Link homology and Frobenius extensions, 
math.QA/0411447. 

\noindent
[KR] M.~Khovanov and L.~Rozansky, Matrix factorizations and
link homology, math.QA 0401268. 

\noindent
[L] E.~S.~Lee, An endomorphism of the Khovanov invariant, 
 math.GT/0210213. 

\noindent 
[OS] P.~Ozsv\'ath and Z.~Szab\'o, Holomorphic disks and 
knot invariants, arXiv math.GT/0209056. 

\noindent
[PT] J.~Przytycki and P.~Traczyk, Conway Algebras and 
Skein Equivalence of Links, \emph{ Proc. AMS}  100, 744-748, 1987.

\noindent
[R1] J.~Rasmussen, Floer homology and knot complements, 
math.GT/0306378. 

\noindent
[R2] J.~Rasmussen, Khovanov homology and the slice genus, 
math.GT/0402131. 

\noindent 
[V] O.~Viro, Private communication. 

\vspace{0.2in} 

\noindent 
Mikhail Khovanov, Department of Mathematics, University of California,
Davis, CA 95616, mikhail@math.ucdavis.edu
 
\vspace{0.07in}
 
\noindent 
Lev Rozansky, Department of Mathematics, University of North Carolina
Chapel Hill, NC 27599, rozansky@math.unc.edu

\end{document}